\documentclass[oneside,draft,11pt]{amsart}

\newcommand{\R}{\mathds R}
\newcommand{\dd}{\mathrm d}
\newcommand{\bfG}{\mathbf{\Gamma}}
\newcommand{\nablao}{\nabla^{\hskip-0.8pt\scriptscriptstyle\otimes}}
\newcommand{\FR}{\mathrm{FR}}
\newcommand{\FRo}{\mathrm{FR}^{\mathrm o}}
\newcommand{\FRu}{\mathrm{FR}^{\mathrm u}}
\newcommand{\FRc}{\mathrm{FR}^{\mathrm c}}
\newcommand{\Or}{\mathrm O}
\newcommand{\SO}{\mathrm{SO}}
\newcommand{\Ur}{\mathrm U}
\newcommand{\so}{\mathfrak{so}}
\newcommand{\GL}{\mathrm{GL}}
\newcommand{\gl}{\mathfrak{gl}}
\newcommand{\Hor}{\mathrm{Hor}}
\newcommand{\Ad}{\mathrm{Ad}}
\newcommand{\Id}{\mathrm{Id}}
\newcommand{\pip}{\Pi}
\newcommand{\pr}{\mathrm{pr}}
\newcommand{\In}[1]{\mathfrak I^{#1}}
\newcommand{\bfd}{{\mathrm{d\!I}}}
\newcommand{\Sym}{\mathrm{sym}}
\newcommand{\Lin}{\mathrm{Lin}}
\newcommand{\alt}{\mathrm a}
\newcommand{\Dim}{\mathrm{dim}}
\newcommand{\Tor}{\mathrm T}
\newcommand{\rank}{\mathrm{rank}}
\newcommand{\oo}{\mathrm o}

\usepackage{times}
\usepackage{dsfont}
\usepackage[all]{xy}

\title[$G$-structure preserving affine immersions]{An existence theorem for $\mathbf G$-structure preserving affine immersions}

\author[P. Piccione]{Paolo Piccione}
\author[D. Tausk]{Daniel V. Tausk}

\address{Departamento de Matem\'atica,\hfill\break\indent  Universidade de S\~ao Paulo, Brazil}
\email{piccione@ime.usp.br, tausk@ime.usp.br}
\urladdr{http://www.ime.usp.br/\~{}piccione,\hfill\break\phantom{URL: } http://www.ime.usp.br/\~{}tausk}

\subjclass[2000]{53A15, 53B05, 53C10, 53C40}

\keywords{$G$-structures, affine immersions}

\date{October 2007}

\begin{document}

\makeatletter
\renewenvironment{proof}[1][\proofname]{\par
  \pushQED{\qed}%
  \normalfont \topsep6\p@\@plus6\p@\relax
  \trivlist
  \item[\hskip\labelsep
        \scshape
    #1\@addpunct{.}]\ignorespaces
}{%
  \popQED\endtrivlist\@endpefalse
}
\makeatother

\numberwithin{equation}{section}
\theoremstyle{plain}\newtheorem{teo}{Theorem}[section]
\theoremstyle{plain}\newtheorem{prop}[teo]{Proposition}
\theoremstyle{plain}\newtheorem{lem}[teo]{Lemma}
\theoremstyle{plain}\newtheorem{cor}[teo]{Corollary}
\theoremstyle{definition}\newtheorem{defin}[teo]{Definition}
\theoremstyle{remark}\newtheorem{rem}[teo]{Remark}
\theoremstyle{plain} \newtheorem{assum}[teo]{Assumption}
\swapnumbers
\theoremstyle{definition}\newtheorem{example}{Example}[section]

\begin{abstract}
We prove an existence result for local and global $G$-structure preserving affine immersions between
affine manifolds. Several examples are discussed in the context of Riemannian and semi-Riemannian
geometry, including the case of isometric immersions into Lie groups endowed with
a left-invariant metric, and the case of isometric immersions into products of space forms.
\end{abstract}

\maketitle

\begin{section}{Introduction}
In this paper we prove an existence result for $G$-structure preserving affine immersions between
affine manifolds, with special attention to the class of isometric immersions in the context
of Riemannian and semi-Riemannian geometry.
The original idea was to find a unifying language for several isometric
immersion theorems that appear in the classical literature \cite{Dajczer}
(immersions into Riemannian manifolds with constant sectional curvature,
immersions into K\"ahler manifolds of constant holomorphic curvature),
and also some recent results (see for instance \cite{Benoit1, Benoit}) concerning
the existence of isometric immersions in more general Riemannian
manifolds. Given an isometric immersion, the celebrated equations of Gauss, Codazzi and Ricci relate
the curvature of the environment with the curvature of the submanifold, the curvature of the normal bundle
and the second fundamental form (and its covariant derivative). A folk theorem says that such equations are
necessary conditions for the existence of isometric immersions, however the reader should observe that,
unless the isometric immersion has already been given, the equations cannot in general even be written down.
Nevertheless, when the target manifold is ``sufficiently homogeneous'' (or, more precisely,
{\em infinitesimally homogeneous\/} in the sense of Definition~\ref{thm:definfhom}), the Gauss, Codazzi
and Ricci equations do make sense {\it a priori\/} and then they are indeed necessary conditions for the existence
of the isometric immersion. Sufficient conditions for the existence of an isometric immersion involve
additional assumptions depending on the context; the starting point
of our theory was precisely the interpretation of such additional
assumptions in terms of $G$-structures and {\em inner torsion}, which is a kind of covariant derivative of
a $G$-structure.

The central result of the paper is an affine immersion theorem into infinitesimally
homogeneous affine manifolds endowed with a $G$-structure. {\em Infinitesimally homogeneous\/}
means that the curvature and the torsion of the connection, as well
as the inner torsion of the $G$-structure, can be written
uniquely in terms of the $G$-structure, i.e., are constant in frames that belong to
the $G$-structure. For instance, consider the case that $M$ is a Riemannian
manifold endowed with the Levi-Civita connection of its metric tensor, $G$
is the orthogonal group and the $G$-structure is given by the set of orthonormal
frames. Since the Levi-Civita connection is compatible with the Riemannian metric,
the inner torsion of this $G$-structure is zero. The condition that the
curvature tensor should be constant in orthonormal frames is equivalent
to the condition that $M$ has constant sectional curvature, and we recover
in this case the classical ``fundamental theorem of isometric immersions
in spaces of constant curvature''. Similarly, if $M$ is a Riemannian manifold
endowed with an orthogonal almost complex structure, then one has a $G$-structure
on $M$, where $G$ is the unitary group, by considering the set of orthonormal
complex frames of $TM$. In this case, the inner torsion of the $G$-structure
relatively to the Levi-Civita connection of the Riemannian metric is
the covariant derivative of the almost complex structure, which vanishes if and only
if $M$ is K\"ahler. Requiring that the
curvature tensor be constant in orthonormal complex frames means that
$M$ has constant holomorphic curvature; in this context, our immersion theorem
reproduces the classical result of isometric immersions into K\"ahler manifolds of
constant holomorphic curvature. Another interesting example of $G$-structure
that will be considered in some detail is the case of Riemannian manifolds
endowed with a distinguished unit vector field $\xi$; in this case,
we obtain an immersion theorem into Riemannian manifolds with the property that both
the curvature tensor and the covariant derivative of the vector field at a general point $p$
can be written in terms only of the Riemannian metric at $p$ and of the vector $\xi(p)$. This is the case in a number
of important examples, like for instance all manifolds that are Riemannian products
of a space form with a copy of the real line, as well as all homogeneous, simply-connected
$3$-dimensional manifolds whose isometry group has dimension $4$. These examples were first considered
in \cite{Benoit1}. Two more examples will be studied in some detail. First, we
will consider isometric immersions into Lie groups endowed with a left
invariant semi-Riemannian metric tensor. These manifolds have an obvious $1$-structure,
given by the choice of a distinguished orthonormal left invariant frame; clearly, the curvature
tensor is constant in this frame. Moreover, the inner torsion of the structure is
simply the Christoffel tensor associated to this frame, which is also constant.
It should be observed that a different immersion theorem into a class of nilpotent and solvable Lie groups has been
recently proved in \cite{JorgeHerbert}. In spite of many analogies both in the statement
and in the proof of the result, the setup considered by the author in \cite{JorgeHerbert} does not
fit into the infinitesimally homogeneous case considered in the present paper.
Another example discussed is the case of isometric immersions
into products of manifolds with constant sectional curvature; in this situation, the $G$-structure
considered is the one consisting of orthonormal frames adapted to such product. More generally,
products of infinitesimally homogeneous affine manifolds with $G$-structures are infinitesimally homogeneous.

The proof of (the local version of) the main theorem relies on an application of Frobenius in the language
of differential forms. More precisely, assume that we are given affine manifolds $(M,\nabla)$, $(\overline M,\overline\nabla)$
and a vector bundle $E^0$ over $M$. Given a suitable set of data (a connection $\nabla^0$ on $E^0$,
and ``second fundamental forms'' $\alpha^0$ and $A^0$) we assemble a connection $\widehat\nabla$ on
the Whitney sum $\widehat E=TM\oplus E$ and we look for an immersion $f:M\to\overline M$ and a connection
preserving vector bundle isomorphism $L:\widehat E\to f^*T\overline M$ such that $L\vert_{TM}=\dd f$. We assume that $G$-structures
$\widehat P$ and $\overline P$ are given on $\widehat E$ and on $T\overline M$, respectively, and we require
that $L$ be $G$-structure preserving. Given a smooth local frame of $\widehat E$ (in $\widehat P$), the problem
of determining $L$ is reduced to the problem of determining a smooth map $F:U\to\overline P$, with $U$ open in $M$,
such that $F$ pulls back the canonical form and the connection form of $\overline P$ to, respectively, the canonical
form and the connection form of $\widehat P$. We then employ a version of the Frobenius theorem that allows
one to guarantee the existence of a smooth map $F$ satisfying a PDE of the form $F^*\lambda^{\overline M}=\lambda^M$, where
$\lambda^M$, $\lambda^{\overline M}$ are vector-valued $1$-forms taking values in the same vector space. The integrability
condition for such PDE corresponds to the Gauss, Codazzi and Ricci equations, as well as to some {\em torsion equations};
moreover, the condition that $\lambda^M$, $\lambda^{\overline M}$ take value in the same vector space corresponds
to an equation relating the inner torsions of $\widehat P$ and of $\overline P$.
Finally, the proof of the global version of the affine immersion theorem employs a general globalization principle stated in the language
of pre-sheafs.
\end{section}

\begin{section}{Notations and terminology}

\subsection*{Vector spaces}
Let $V$ be a real finite-dimensional vector space. We denote by $\GL(V)$ the general linear group of $V$ and by
$\gl(V)$ its Lie algebra. If $W$ is another real finite-dimensional vector space and $p:V\to W$ is a linear isomorphism
then $\mathcal I_p:\GL(V)\to\GL(W)$ denotes the Lie group isomorphism given by conjugation with $p$ and
$\Ad_p=\dd\mathcal I_p(\Id):\gl(V)\to\gl(W)$ denotes the Lie algebra isomorphism given by conjugation with $p$.
By $\Lin(V,W)$ we denote the space of linear maps from $V$ to $W$.

\subsection*{Vector bundles, frame bundles and connections}
Let $E$ be a vector bundle over a differentiable manifold $M$. We denote by $\bfG(E)$ the set of all smooth sections of $E$.
Given a connection $\nabla$ on $E$ then the {\em curvature\/} of $\nabla$ is the smooth tensor $R\in\bfG(TM^*\otimes TM^*\otimes E^*\otimes E)$
defined by:
\[R(X,Y)\epsilon=\nabla_X\nabla_Y\epsilon-\nabla_Y\nabla_X\epsilon-\nabla_{[X,Y]}\epsilon,\]
for all $X,Y\in\bfG(TM)$, $\epsilon\in\bfG(E)$; if $\iota:TM\to E$ is a vector bundle morphism then the
{\em $\iota$-torsion\/} of $\nabla$ is the smooth tensor $\Tor^\iota\in\bfG(TM^*\otimes TM^*\otimes E)$ defined
by:
\[\Tor^\iota(X,Y)=\nabla_X\big(\iota(Y)\big)-\nabla_Y\big(\iota(X)\big)-\iota\big([X,Y]\big),\]
for all $X,Y\in\bfG(TM)$. When $E=TM$ and $\iota$ is the identity $\Tor=\Tor^\iota$ coincides with the usual {\em torsion\/}
of $\nabla$.

Let $k$ be the {\em rank\/} of $E$, i.e., the dimension of the fibers of $E$. We denote by $\FR(E)=\bigcup_{x\in M}\FR(E_x)$
the frame bundle of $E$, which is the set of all linear isomorphisms $p:\R^k\to E_x$, with $x\in M$.
The frame bundle $\FR(E)$ is a $\GL(\R^k)$-principal bundle over $M$. A local section $s:U\to\FR(E)$ (where
$U$ is an open subset of $M$) is called a {\em local frame\/} for $E$. A smooth local frame $s:U\to\FR(E)$
defines a connection $\bfd^s$ on $E\vert_U$ which corresponds via the trivialization of $E\vert_U$ defined
by $s$ to the standard derivative. More explicitly, we set:
\[\bfd^s_v\epsilon=s(x)\big(\dd\tilde\epsilon_x(v)\big),\]
for all $x\in U$, $v\in T_xM$, $\epsilon\in\bfG(E\vert_U)$, where $\tilde\epsilon:U\to\R^k$ is defined
by $\tilde\epsilon(x)=s(x)^{-1}\big(\epsilon(x)\big)$, for all $x\in U$. If $\nabla$ is a connection in $E$
then the {\em Christoffel tensor\/} of $\nabla$ with respect to the smooth local frame $s$ is the tensor
$\Gamma=\nabla-\bfd^s$; more explicitly, $\Gamma:U\to TM^*\otimes E^*\otimes E$ is the smooth local section
such that:
\[\nabla_v\epsilon=\bfd^s_v\epsilon+\Gamma_x\big(v,\epsilon(x)\big),\]
for all $x\in U$, $v\in T_xM$ and all $\epsilon\in\bfG(E\vert_U)$. For $x\in U$, $v\in T_xM$, we also
set $\Gamma_x(v)=\Gamma_x(v,\cdot)\in\gl(E_x)$, so that $\Gamma_x:T_xM\to\gl(E_x)$ is a linear map.
If $\Hor$ is the horizontal distribution in $\FR(E)$ corresponding
to $\nabla$ and $\omega$ is the $\gl(\R^k)$-valued connection form on $\FR(E)$ whose kernel is $\Hor$ then,
setting $\bar\omega=s^*\omega$, we have:
\begin{equation}\label{eq:relbaromegaGamma}
\Gamma_x(v)=\Ad_{s(x)}\big(\bar\omega_x(v)\big),
\end{equation}
for all $x\in U$ and all $v\in T_xM$. If $\iota:TM\to E$ is a vector bundle morphism then the {\em $\iota$-canonical
form\/} of $\FR(E)$ is the $\R^k$-valued $1$-form $\theta$ on $\FR(E)$ defined by:
\[\theta_p(\zeta)=p^{-1}\big[\iota\big(\dd\pip_p(\zeta)\big)\big]\in\R^k,\]
for all $p\in\FR(E)$, $\zeta\in T_p\FR(E)$, where $\pip:\FR(E)\to M$ denotes the projection. When $E=TM$
and $\iota$ is the identity then $\theta$ is simply the usual {\em canonical form\/} of $\FR(TM)$. The
{\em $\iota$-torsion form\/} $\Theta$ and the {\em curvature form\/} $\Omega$ are defined respectively by:
\begin{equation}\label{eq:defThetaOmega}
\Theta=\dd\theta+\omega\wedge\theta,\quad\Omega=\dd\omega+\omega\wedge\omega,
\end{equation}
where the wedge product in $\omega\wedge\omega$ is taken with respect to the associative product of $\gl(\R^k)$.
The following equalities hold:
\begin{gather}
\label{eq:ThetaT}p\big(\Theta_p(\zeta_1,\zeta_2)\big)=\Tor^\iota_x\big(\dd\pip_p(\zeta_1),\dd\pip_p(\zeta_2)\big)\in E_x,\\
\label{eq:OmegaR}\Ad_p\big(\Omega_p(\zeta_1,\zeta_2)\big)=R_x\big(\dd\pip_p(\zeta_1),\dd\pip_p(\zeta_2)\big)\in\gl(E_x),
\end{gather}
for all $p\in\FR(E)$, $\zeta_1,\zeta_2\in T_p\FR(E)$, where $x=\pip(p)$.

\subsection*{Covariant derivative along curves}
Let $p:I\to\FR(E)$ be a smooth curve and set $\gamma=\pip\circ p$, where $\pip:\FR(E)\to M$ denotes the projection.
For all $t\in I$, we denote by $(\nabla_1p)(t)\in\Lin(\R^k,E_{\gamma(t)})$ the {\em covariant
derivative\/} of $p$ at the instant $t$, which is just the vertical component of $p'(t)\in T_{p(t)}\FR(E)$ (observe
that the vertical space of $\FR(E)$ at $p(t)$ is identified with $\Lin(\R^k,E_{\gamma(t)})$). We have:
\begin{equation}\label{eq:omegaplinha}
\omega_{p(t)}\big(p'(t)\big)=p(t)^{-1}\circ(\nabla_1p)(t),
\end{equation}
for all $t\in I$. If $\varepsilon:I\to E$ is a smooth section of $E$ along $\gamma$
(i.e., $\varepsilon(t)\in E_{\gamma(t)}$, for all $t\in I$),
we denote by $(\nabla_1\varepsilon)(t)\in E_{\gamma(t)}$ the covariant derivative of $\varepsilon$
along $\gamma$ at the instant $t$. Given a smooth curve $u:I\to\R^k$ then $\varepsilon:I\ni t\mapsto p(t)\cdot u(t)\in E_{\gamma(t)}$
is a smooth section of $E$ along $\gamma$ and the following ``Leibniz rule'' holds:
\[(\nabla_1\varepsilon)(t)=(\nabla_1p)(t)\cdot u(t)+p(t)\cdot u'(t),\]
for all $t\in I$. If $s:U\to\FR(E)$ is a smooth local section with $\gamma(I)\subset U$ and if
$\bar\omega=s^*\omega$ then:
\begin{equation}\label{eq:nabla1p}
(\nabla_1p)(t)=s\big(\gamma(t)\big)\circ\big[\tilde p'(t)+\bar\omega_{\gamma(t)}\big(\gamma'(t)\big)\circ\tilde p(t)\big],
\end{equation}
for all $t\in I$, where $\tilde p:I\to\GL(\R^k)$ is defined by $\tilde p(t)=s\big(\gamma(t)\big)^{-1}\circ p(t)$.

\subsection*{Vector subbundles}
If $F$ is a vector subbundle of $E$ then the {\em absolute second fundamental form\/} of $F$ in $E$ with
respect to the connection $\nabla$ is the tensor $\alpha^F\in\bfG\big(TM^*\otimes F^*\otimes(E/F)\big)$
defined by $\alpha^F(X,\epsilon)=\mathfrak q\big(\nabla_X\epsilon)$, for all $X\in\bfG(TM)$, $\epsilon\in\bfG(F)$,
where $\mathfrak q:E\to E/F$ denotes the quotient map.

\end{section}

\begin{section}{Affine immersions and their invariants}

Let $(M,\nabla)$, $(\overline M,\overline\nabla)$ be affine manifolds and $f:M\to\overline M$ be a smooth immersion.
We identify the differential $\dd f:TM\to T\overline M$ with a (injective) vector bundle morphism
$\dd f:TM\to f^*T\overline M$.
A vector subbundle $E$ of $f^*T\overline M$ with $f^*T\overline M=\dd f(TM)\oplus E$ will be called a {\em normal
bundle\/} for $f$. Let a normal bundle $E$ for $f$ be fixed; we denote by $\pi^E:f^*T\overline M\to E$
the projection onto $E$ corresponding to the decomposition $f^*T\overline M=\dd f(TM)\oplus E$ and by
$\pi^{TM}:f^*T\overline M\to TM$ the composition of the projection onto $\dd f(TM)$ with
the isomorphism $\dd f^{-1}:\dd f(TM)\to TM$.
Given smooth vector fields $X,Y\in\bfG(TM)$ in $M$, we set $\alpha(X,Y)=\pi^E\big(\overline\nabla_X\,\dd f(Y)\big)\in\bfG(E)$,
so that $\alpha$ is identified with a smooth section of $TM^*\otimes TM^*\otimes E$. We call $\alpha$ the
{\em second fundamental form\/} of the immersion $f$ with respect to the normal bundle $E$.
Notice that in the case of Riemannian (or semi-Riemannian) geometry, $f$ has a canonical normal bundle
(the orthogonal complement of $\dd f(TM)$ with respect to the metric),
so there is also a canonical notion of second fundamental form. Also, if $(M,g)$ and $(\overline M,\bar g)$
are semi-Riemannian, $\nabla$ and $\overline\nabla$ are the corresponding Levi-Civita connections,
$f$ is an isometric immersion and $E$ is the orthogonal complement of $\dd f(TM)$ in $f^*T\overline M$ with respect to $\bar g$ then:
\begin{equation}\label{eq:nablaalpha}
\overline\nabla_X\,\dd f(Y)=\dd f(\nabla_XY)+\alpha(X,Y),
\end{equation}
for all $X,Y\in\bfG(TM)$. In the general affine case, we say that $f$ is an {\em affine immersion with respect
to $E$\/} if \eqref{eq:nablaalpha} holds, for all $X,Y\in\bfG(TM)$. Following \cite{NomizuSasaki}
we say simply that $f$ is an {\em affine immersion\/} if there exists a normal bundle $E$ for $f$
such that $f$ is an affine immersion with respect to $E$. We define the {\em normal connection\/} $\nabla^\perp$
of the immersion $f$ corresponding to the normal bundle $E$ by setting
$\nabla^\perp_X\epsilon=\pi^E(\overline\nabla_X\epsilon)$, for all $X\in\bfG(TM)$ and all
$\epsilon\in\bfG(E)$. In the semi-Riemannian case, $\alpha$ and $\nabla^\perp$ are the only invariants associated to an
isometric immersion. In the general affine case, we have an additional invariant associated to the immersion. We set
$A(X,\epsilon)=\pi^{TM}(\overline\nabla_X\epsilon)$, for all $X\in\bfG(TM)$ and all $\epsilon\in\bfG(E)$;
clearly $A$ is identified with a smooth section of $TM^*\otimes E^*\otimes TM$. We call $A$
the {\em Weingarten form\/} of the immersion $f$ relatively to $E$ and for all $x\in M$ and all $e\in E_x$,
the linear endomorphism $A(e)=A_x(\cdot,e):T_xM\to T_xM$ is called the {\em Weingarten operator\/} in the direction
of $e$. In the semi-Riemannian case $\alpha$ and $A$ are related by the equality:
\[\bar g_x\big(\alpha_x(v,w),e\big)=-\bar g_x\big(A(e)\cdot v,w\big),\quad x\in M,\ v,w\in T_xM,\ e\in E_x,\]
so that $A$ is determined by $\alpha$. In the affine case, there is no relation between $\alpha$ and $A$; moreover,
$\alpha$ is not in general symmetric unless the connection $\overline\nabla$ is symmetric.

We are interested in studying the existence of affine immersions with prescribed invariants $\nabla^\perp$, $\alpha$
and $A$. More precisely, let $(M,\nabla)$, $(\overline M,\overline\nabla)$ be affine manifolds, $E^0$ be a vector bundle over $M$,
$\nabla^0$ be a connection in $E^0$ and $\alpha^0$, $A^0$ be smooth sections of $TM^*\otimes TM^*\otimes E^0$
and $TM^*\otimes(E^0)^*\otimes TM$ respectively. We look for an affine immersion $f:M\to\overline M$,
a normal bundle $E$ for $f$ and a connection preserving vector bundle isomorphism $S:(E^0,\nabla^0)\to(E,\nabla^\perp)$ such that
$S\big(\alpha^0(\cdot,\cdot)\big)=\alpha$ and $A(\cdot,S\cdot)=A^0$. The pair $(f,S)$ will be called a {\em solution
for the affine immersion problem with data $\nabla^0$, $\alpha^0$ and $A^0$}. More generally,
if $f:U\to\overline M$ is an affine immersion defined in an open subset $U$ of $M$ and
$S:E^0\vert_U\to E\vert_U$ is a connection preserving vector bundle isomorphism such that
$S\big(\alpha^0(\cdot,\cdot)\big)=\alpha$ and $A(\cdot,S\cdot)=A^0$ then the pair $(f,S)$ will be called a
{\em local solution for the affine immersion problem with data $\nabla^0$, $\alpha^0$, $A^0$ and with domain $U$}.

An important special situation is the one of isometric immersions. Assume
that $(M,g)$, $(\overline M,\bar g)$ are semi-Riemannian manifolds,
$E^0$ is a vector bundle over $M$ endowed with a semi-Riemannian
structure $g^0$ (i.e., $g^0$ is a smooth section of $(E^0)^*\otimes(E^0)^*$ and $g^0_x$ is a nondegenerate
symmetric bilinear form on $E^0_x$, for all $x\in M$), $g^0$ is $\nabla^0$-parallel and $\alpha^0$
is a smooth symmetric section of the vector bundle $TM^*\otimes TM^*\otimes E^0$; by a {\em solution for the isometric
immersion problem with data $\nabla^0$, $\alpha^0$, $g^0$\/} we mean a pair $(f,S)$ where
$f:M\to\overline M$ is an isometric immersion, $S:(E^0,\nabla^0,g^0)\to(E,\nabla^\perp,g^\perp)$ is a connection preserving
vector bundle isometry and $S\big(\alpha^0(\cdot,\cdot)\big)=\alpha$, where $E$ denotes the orthogonal
complement of $\dd f(TM)$ in $f^*T\overline M$ with respect to $\bar g$ and $g^\perp$ denotes the restriction of $\bar g$
to $E$. As in the affine case, one defines the concept of {\em local\/} solution for the isometric immersion problem by replacing $M$
with an open subset $U$ of $M$.

Notice that if $(f,S)$ is a (local) solution for the isometric immersion problem, $(M,g)$
and $(\overline M,\bar g)$ are endowed with their respective Levi-Civita connections $\nabla$, $\overline\nabla$
and a smooth section $A^0$ of $TM^*\otimes(E^0)^*\otimes TM$ is defined by the equality:
\begin{equation}\label{eq:g0alpha0A0}
g^0_x\big(\alpha^0_x(v,w),e\big)=-g_x\big(A^0_x(e)\cdot v,w\big),\quad x\in M,\ v,w\in T_xM,\ e\in E^0_x,
\end{equation}
then $(f,S)$ is also a (local) solution for the affine immersion problem with data $\nabla^0$, $\alpha^0$ and $A^0$.

\end{section}

\begin{section}{The components of a connection}

The following general construction gives a convenient language for discussing the theory of affine immersions.
Let $\mathcal E$ be a vector bundle over a differentiable manifold $M$ endowed with a connection $\nabla$ and
let $\mathcal E=\mathcal E^1\oplus\mathcal E^2$ be a direct sum decomposition of $\mathcal E$; denote
by $\pi^i:\mathcal E\to\mathcal E^i$, $i=1,2$, the projections. We set:
\begin{gather*}
\nabla^i_X\epsilon_i=\pi^i(\nabla_X\epsilon_i),\quad i=1,2\\
\alpha^2(X,\epsilon_1)=\pi^2(\nabla_X\epsilon_1),\quad\alpha^1(X,\epsilon_2)=\pi^1(\nabla_X\epsilon_2),
\end{gather*}
for all $X\in\bfG(TM)$, $\epsilon_i\in\bfG(\mathcal E^i)$, $i=1,2$, so that $\nabla^i$ is a connection in $\mathcal E^i$,
$i=1,2$, and $\alpha^2$ (resp., $\alpha^1$) is identified with a smooth section of $TM^*\otimes(\mathcal E^1)^*\otimes\mathcal E^2$
(resp., of $TM^*\otimes(\mathcal E^2)^*\otimes\mathcal E^1$). We call $\nabla^1$, $\nabla^2$, $\alpha^2$ and $\alpha^1$
the {\em components\/} of the connection $\nabla$ with respect to the decomposition $\mathcal E=\mathcal E^1\oplus\mathcal E^2$.
Clearly, one recovers $\nabla$ from its components using the formula:
\[\nabla_X\epsilon=\nabla^1_X\big(\pi^1(\epsilon)\big)+\alpha^1\big(X,\pi^2(\epsilon)\big)
+\nabla^2_X\big(\pi^2(\epsilon)\big)+\alpha^2\big(X,\pi^1(\epsilon)\big),\]
where $X\in\bfG(TM)$ and $\epsilon\in\bfG(\mathcal E)$.

Notice that if $f:(M,\nabla)\to(\overline M,\overline\nabla)$ is an affine immersion with respect
to a normal bundle $E$ then, identifying for a moment $TM$ with $\dd f(TM)$,
the components of the connection $f^*\overline\nabla$ on $f^*T\overline M$ with
respect to the decomposition $f^*T\overline M=\dd f(TM)\oplus E$ are the connection $\nabla$ of $M$, the normal connection
$\nabla^\perp$, the second fundamental form $\alpha$ and the Weingarten form $A$ of $f$ with respect to $E$.

\begin{rem}\label{thm:remonly3components}
If $g^1$, $g^2$ are semi-Riemannian structures on $E^1$ and $E^2$, respectively and if $g$ is the semi-Riemannian
structure on $E$ given by the orthogonal direct sum of $g^1$ and $g^2$ then a connection $\nabla$ with
components $\nabla^1$, $\nabla^2$, $\alpha^2$, $\alpha^1$ is compatible with $g$ (i.e., $\nabla g=0$) if and only
if $\nabla^i$ is compatible with $g^i$, $i=1,2$, and the following relation between $\alpha^2$ and $\alpha^1$ holds:
\begin{equation}\label{eq:relationalphas}
g_x\big(\alpha^2_x(v,e_1),e_2\big)+g_x\big(e_1,\alpha^1_x(v,e_2)\big)=0,
\end{equation}
for all $x\in M$, $v\in T_xM$, $e_1\in\mathcal E^1_x$ and $e_2\in\mathcal E^2_x$. Notice that relation
\eqref{eq:relationalphas} implies that $\alpha^1$ is uniquely determined from $\alpha^2$, so that in a context
where we are dealing with connections compatible with a semi-Riemannian structure, we will talk only about the components
$\nabla^1$, $\nabla^2$ and $\alpha^2$ of $\nabla$, where one should understand implicitly that $\alpha^1$ is determined
by condition \eqref{eq:relationalphas}.
\end{rem}

Denote by $R$, $R^1$, $R^2$ the curvature tensors of $\nabla$, $\nabla^1$ and $\nabla^2$, respectively.
A straightforward computation gives the following:
\begin{gather}
\label{eq:genGauss}\pi^1\big(R_x(v,w)e_1\big)=R^1_x(v,w)e_1+\alpha^1_x\big(v,\alpha^2_x(w,e_1)\big)
-\alpha^1_x\big(w,\alpha^2_x(v,e_1)\big),\\
\label{eq:genRicci}\pi^2\big(R_x(v,w)e_2\big)=R^2_x(v,w)e_2+\alpha^2_x\big(v,\alpha^1_x(w,e_2)\big)
-\alpha^2_x\big(w,\alpha^1_x(v,e_2)\big),
\end{gather}
for all $x\in M$, $e_1\in\mathcal E^1_x$, $e_2\in\mathcal E^2_x$ and all $v,w\in T_xM$. Moreover,
given a connection $\nabla^M$ on $TM$ with torsion $\Tor$ and denoting by $\nablao$
the induced connections on $TM^*\otimes(\mathcal E^2)^*\otimes\mathcal E^1$ and on $TM^*\otimes(\mathcal E^1)^*\otimes\mathcal E^2$ then:
\begin{gather}
\label{eq:genCodazzi}\begin{aligned}
\pi^2\big(R_x(v,w)e_1\big)=(\nablao\alpha^2)_x(v,w,e_1)
&-(\nablao\alpha^2)_x(w,v,e_1)\\
&+\alpha^2_x\big(\Tor_x(v,w),e_1\big),
\end{aligned}\\[5pt]
\label{eq:genCodazzi2}\begin{aligned}
\pi^1\big(R_x(v,w)e_2\big)=(\nablao\alpha^1)_x(v,w,e_2)
&-(\nablao\alpha^1)_x(w,v,e_2)\\
&+\alpha^1_x\big(\Tor_x(v,w),e_2\big),
\end{aligned}
\end{gather}
for all $x\in M$, $e_1\in\mathcal E^1_x$, $e_2\in\mathcal E^2_x$ and all $v,w\in T_xM$. If
$\iota=(\iota^1,\iota^2):TM\to\mathcal E=\mathcal E^1\oplus\mathcal E^2$ is a vector bundle morphism
then the $\iota$-torsion $\Tor^\iota$ of $\nabla$ satisfies the following identities:
\begin{gather}
\label{eq:gentorsioneq1}\pi^1\big(\Tor^\iota_x(v,w)\big)=\Tor^{\iota^1}_x(v,w)+\alpha^1_x\big(v,\iota^2(w)\big)
-\alpha^1_x\big(w,\iota^2(v)\big),\\
\label{eq:gentorsioneq2}\pi^2\big(\Tor^\iota_x(v,w)\big)=\Tor^{\iota^2}_x(v,w)+\alpha^2_x\big(v,\iota^1(w)\big)
-\alpha^2_x\big(w,\iota^1(v)\big),
\end{gather}
for all $x\in M$, $v,w\in T_xM$, where $\Tor^{\iota^1}$, $\Tor^{\iota^2}$ denote respectively the $\iota^1$-torsion
of $\nabla^1$ and the $\iota^2$-torsion of $\nabla^2$.

Let us look at equations \eqref{eq:genGauss}, \eqref{eq:genRicci}, \eqref{eq:genCodazzi}, \eqref{eq:genCodazzi2},
\eqref{eq:gentorsioneq1} and \eqref{eq:gentorsioneq2} in the context of affine immersions.
More precisely, let $f:(M,\nabla)\to(\overline M,\overline\nabla)$
be an affine immersion, $E$ be a normal bundle for $f$, and $\nabla^\perp$, $\alpha$ and $A$ denote respectively
the normal connection, the second fundamental form and the Weingarten form. Since $\nabla$, $\nabla^\perp$,
$\alpha$ and $A$ are (up to the identification of $TM$ with $\dd f(TM)$)
the components of $f^*\overline\nabla$ with respect to the decomposition
$f^*T\overline M=\dd f(TM)\oplus E$, equation \eqref{eq:genGauss} gives:
\begin{multline}\label{eq:Gauss}
\pi^{TM}\big[\overline R_{f(x)}\big(\dd f_x(v),\dd f_x(w)\big)\dd f_x(u)\big]=R_x(v,w)u+A_x\big(v,\alpha_x(w,u)\big)\\
-A_x\big(w,\alpha_x(v,u)\big),
\end{multline}
for all $x\in M$ and all $v,w,u\in T_xM$, where $R$ and $\overline R$ denote the curvature tensors
of $\nabla$ and $\overline\nabla$ respectively. We call \eqref{eq:Gauss} the {\em Gauss equation\/}
of the affine immersion $f$ with respect to $E$. Similarly, equation \eqref{eq:genRicci} gives:
\begin{multline}\label{eq:Ricci}
\pi^E\big[\overline R_{f(x)}\big(\dd f_x(v),\dd f_x(w)\big)e\big]=R^\perp_x(v,w)e+\alpha_x\big(v,A_x(w,e)\big)\\
-\alpha_x\big(w,A_x(v,e)\big),
\end{multline}
for all $x\in M$, $v,w\in T_xM$ and all $e\in E_x$, where $R^\perp$ denotes the curvature tensor of
the normal connection $\nabla^\perp$. We call \eqref{eq:Ricci} the {\em Ricci equation\/}
of the affine immersion $f$ with respect to $E$. Equations \eqref{eq:genCodazzi}
and \eqref{eq:genCodazzi2} (with $\nabla^M=\nabla$) give:
\begin{gather}
\label{eq:Codazzi}\begin{aligned}
\pi^E\big[\overline R_{f(x)}\big(\dd f_x(v),\dd f_x(w)\big)u\big]=(\nablao\alpha)_x(v,w,u)
&-(\nablao\alpha)_x(w,v,u)\\
&+\alpha_x\big(\Tor_x(v,w),u\big),
\end{aligned}\\[5pt]
\label{eq:Codazzi2}\begin{aligned}
\pi^{TM}\big[\overline R_{f(x)}\big(\dd f_x(v),\dd f_x(w)\big)e\big]=(\nablao A)_x(v,w,e)
&-(\nablao A)_x(w,v,e)\\
&+A_x\big(\Tor_x(v,w),e\big),
\end{aligned}
\end{gather}
for all $x\in M$, $v,w,u\in T_xM$ and all $e\in E_x$, where $\Tor$ denotes the torsion tensor of $\nabla$.
We call \eqref{eq:Codazzi} and \eqref{eq:Codazzi2} the {\em Codazzi equations\/} of the affine immersion
$f$ with respect to $E$. Finally, if $\iota:TM\to f^*T\overline M=\dd f(TM)\oplus E$ is the map identified with
$\dd f:TM\to T\overline M$ then $\iota^1=\dd f:TM\to\dd f(TM)$, $\iota^2=0$ and
equations \eqref{eq:gentorsioneq1} and \eqref{eq:gentorsioneq2} read\footnote{%
By taking such map $\iota$, the $\iota$-torsion of $f^*\overline\nabla$ is just the pull-back by $f$ of the torsion
of $\overline\nabla$.}:
\begin{gather}
\label{eq:torsioneq1}\pi^{TM}\big[\overline\Tor_{f(x)}\big(\dd f_x(v),\dd f_x(w)\big)\big]=\Tor_x(v,w),\\
\label{eq:torsioneq2}\pi^E\big[\overline\Tor_{f(x)}\big(\dd f_x(v),\dd f_x(w)\big)\big]=\alpha_x(v,w)-\alpha_x(w,v),
\end{gather}
for all $x\in M$, $v,w\in T_xM$, where $\overline\Tor$ denotes the torsion of $\overline\nabla$.
We call \eqref{eq:torsioneq1} and \eqref{eq:torsioneq2} the {\em torsion equations\/} of the affine immersion
$f$ with respect to $E$.

Now assume that we are given affine manifolds $(M,\nabla)$, $(\overline M,\overline\nabla)$, a vector bundle
$E^0$ over $M$ endowed with a connection $\nabla^0$ and smooth sections $\alpha^0$, $A^0$ of the vector bundles
$TM^*\otimes TM^*\otimes E^0$ and $TM^*\otimes(E^0)^*\otimes TM$ respectively. Assume that there exists a local solution $(f,S)$ for the affine
immersion problem with data $\nabla^0$, $\alpha^0$, $A^0$ defined on an open subset $U$ of $M$, where $S:E^0\vert_U\to E\vert_U$
and $E$ is a normal bundle for $f$. Clearly, \eqref{eq:Gauss}, \eqref{eq:Ricci}, \eqref{eq:Codazzi},
\eqref{eq:Codazzi2}, \eqref{eq:torsioneq1} and \eqref{eq:torsioneq2} imply:
\begin{gather}
\label{eq:Gauss0}
\begin{aligned}
\pi^{TM}\big[\overline R_{f(x)}\big(\dd f_x(v),\dd f_x(w)\big)\dd f_x(u)\big]=R_x(v,w)u&+A^0_x\big(v,\alpha^0_x(w,u)\big)\\
&-A^0_x\big(w,\alpha^0_x(v,u)\big),
\end{aligned}\\[5pt]
\label{eq:Ricci0}
\begin{aligned}
\pi^E\big[\overline R_{f(x)}\big(\dd f_x(v),\dd f_x(w)\big)S_x(e)\big]=S_x\big[R^0_x(v,w)&+\alpha^0_x\big(v,A^0_x(w,e)\big)\\
&-\alpha^0_x\big(w,A^0_x(v,e)\big)\big],
\end{aligned}\\[5pt]
\label{eq:Codazzi0}
\begin{aligned}
\pi^E\big[\overline R_{f(x)}\big(\dd f_x(v),&\dd f_x(w)\big)\dd f_x(u)\big]=S_x\big[(\nablao\alpha^0)_x(v,w,u)\\
&-(\nablao\alpha^0)_x(w,v,u)+\alpha^0_x\big(\Tor_x(v,w),u\big)\big],
\end{aligned}\\[5pt]
\label{eq:Codazzi20}
\begin{aligned}
\pi^{TM}\big[\overline R_{f(x)}\big(\dd f_x(v),&\dd f_x(w)\big)S_x(e)\big]=(\nablao A^0)_x(v,w,e)\\
&-(\nablao A^0)_x(w,v,e)+A^0_x\big(\Tor_x(v,w),e\big),
\end{aligned}\\[5pt]
\label{eq:torsioneq10}\pi^{TM}\big[\overline\Tor_{f(x)}\big(\dd f_x(v),\dd f_x(w)\big)\big]=\Tor_x(v,w),\\
\label{eq:torsioneq20}\pi^E\big[\overline\Tor_{f(x)}\big(\dd f_x(v),\dd f_x(w)\big)\big]=S_x\big(\alpha^0_x(v,w)-\alpha^0_x(w,v)\big),
\end{gather}
for all $x\in U$, $v,w,u\in T_xM$ and all $e\in E^0_x$, where $R^0$ denotes the curvature tensor of $\nabla^0$.

Notice that in the case of isometric immersions, the torsion equation \eqref{eq:torsioneq10} is trivial
and \eqref{eq:torsioneq20} says that $\alpha^0$ is symmetric; moreover, using \eqref{eq:g0alpha0A0}, it can
be seen that the Codazzi equations \eqref{eq:Codazzi0} and \eqref{eq:Codazzi20} are equivalent to each other.

One may think that Gauss, Ricci, Codazzi and the torsion equations are ``necessary conditions''
for the existence of a solution $(f,S)$ of the affine immersion problem, although such statement is
obviously meaningless because one cannot write down equations \eqref{eq:Gauss0}, \eqref{eq:Ricci0},
\eqref{eq:Codazzi0}, \eqref{eq:Codazzi20}, \eqref{eq:torsioneq10} and \eqref{eq:torsioneq20} unless $f$ and $S$ are already given. Notice
that in the special case that $(\overline M,\bar g)$ is a semi-Riemannian manifold with constant
sectional curvature $c\in\R$ and $(f,S)$ is a solution for the isometric immersion problem then the lefthand
side of equations \eqref{eq:Gauss0}, \eqref{eq:Ricci0}, \eqref{eq:Codazzi0} and \eqref{eq:Codazzi20}
can be written only in terms of $c$ and $g$, i.e., without using $f$ and $S$; more explicitly, the lefthand
side of \eqref{eq:Gauss0} is $c\big(g_x(w,u)v-g_x(v,u)w\big)$, while the lefthand sides of
\eqref{eq:Ricci0}, \eqref{eq:Codazzi0} and \eqref{eq:Codazzi20} are zero (the possibility of writing
down Gauss, Ricci and Codazzi equations without using $f$ and $S$ depends on the fact that the curvature tensor
$\overline R$ of the target manifold is {\em constant in orthonormal frames}). Thus, in this case,
the Gauss, Ricci and Codazzi equations are indeed necessary conditions for the existence of
a solution of the isometric immersion problem $(f,S)$. By the celebrated fundamental theorem
of isometric immersions into space forms (see, for instance, \cite{Dajczer, Grif, PaTe, Tenenblat}),
the Gauss, Ricci and Codazzi equations are also sufficient conditions for the existence of local
solutions for the isometric immersion problem (provided that one assumes that $\alpha^0$ is symmetric
and that $g^0$ is $\nabla^0$-parallel).

Using the notion of {\em infinitesimally homogeneous affine manifold with $G$-struc\-ture\/} introduced
in Section~\ref{sec:infhom} we will describe a very general situation in which the lefthand
side of \eqref{eq:Gauss0}, \eqref{eq:Ricci0}, \eqref{eq:Codazzi0}, \eqref{eq:Codazzi20}, \eqref{eq:torsioneq10}
and \eqref{eq:torsioneq20} can be described without explicit use of $f$ and $S$.

\end{section}

\begin{section}{$G$-structures and inner torsion}

Let $E$ be a vector bundle of rank $k$ over a differentiable manifold $M$. If $G$ is a Lie subgroup of $\GL(\R^k)$
then by a {\em $G$-structure\/} on $E$ we mean
a $G$-principal subbundle $P$ of $\FR(E)$. By a $G$-structure on $M$ we mean a $G$-structure on the tangent
bundle of $M$. Let $\nabla$ be a connection in $E$. We denote by $\Hor$ the corresponding horizontal
distribution on $\FR(E)$ and by $\omega$ the $\gl(\R^k)$-valued connection form on $\FR(E)$ whose kernel
is $\Hor$. We say that $\nabla$ is {\em compatible\/} with a $G$-structure $P$ if $\Hor_p\subset T_pP$,
for all $p\in P$, i.e., if parallel transport carries frames in $P$ to frames in $P$. In the general case,
there is a tensor that measures the lack of compatibility of $\nabla$ with $P$ called
the {\em inner torsion\/} of $P$ with respect to $\nabla$, which is defined as follows.

For each $x\in M$, denote by $G_x$ the subgroup of $\GL(E_x)$ consisting of {\em $G$-structure preserving maps},
i.e., maps $\sigma:E_x\to E_x$ such that $\sigma\circ p\in P_x$ for some (and hence for all) $p\in P_x$. Clearly
$G_x=\mathcal I_p(G)$, for all $p\in P_x$, so that $G_x$ is a Lie subgroup of $\GL(E_x)$. We denote by
$\mathfrak g_x\subset\gl(E_x)$ the Lie algebra of $G_x$, so that
$\Ad_p(\mathfrak g)=\mathfrak g_x$, for all $p\in P_x$, where $\mathfrak g\subset\gl(\R^k)$
denotes the Lie algebra of $G$. For each $x\in M$ and each $p\in\FR(E_x)$, we can identify the tangent space
$T_p\FR(E)$ with the direct sum $T_xM\oplus\gl(\R^k)$ via the isomorphism $(\dd\pip_p,\omega_p)$,
where $\pip:\FR(E)\to M$ denotes the projection. For $p\in P$, the subspace $\mathcal V_p=(\dd\pip_p,\omega_p)(T_pP)$
of $T_xM\oplus\gl(\R^k)$ corresponding to $T_pP$ satisfies the conditions $\pr_1(\mathcal V_p)=T_xM$ and
$\mathcal V_p\cap\big(\{0\}\oplus\gl(\R^k)\big)=\{0\}\oplus\mathfrak g$, where $\pr_1:T_xM\oplus\gl(\R^k)\to T_xM$ denotes
the first projection. Thus, there exists a unique linear map $\mathcal L_p:T_xM\to\gl(\R^k)/\mathfrak g$ such that:
\begin{equation}\label{eq:defcalVp}
\mathcal V_p=\big\{(v,X)\in T_xM\oplus\gl(\R^k):\mathcal L_p(v)=X+\mathfrak g\big\}.
\end{equation}
If $s:U\to P$ is a smooth local section with $x\in U$ and $\bar\omega=s^*\omega$ then $\mathcal L_p$ is the composition
of $\bar\omega_x:T_xM\to\gl(\R^k)$ with the quotient map $\gl(\R^k)\to\gl(\R^k)/\mathfrak g$.
It follows from the usual properties of connection forms that, given $p,q\in P_x$, the maps $\mathcal L_p$ and
$\mathcal L_q$ are related by $\mathcal L_q=\overline\Ad_g\circ\mathcal L_p$, where $g\in G$ is such that $p=q\circ g$,
and $\overline\Ad_g:\gl(\R^k)/\mathfrak g\to\gl(\R^k)/\mathfrak g$ is obtained from $\Ad_g:\gl(\R^k)\to\gl(\R^k)$
by passing to the quotient. It follows that the linear map $\In P_x:T_xM\to\gl(E_x)/\mathfrak g_x$
defined by:
\begin{equation}\label{eq:relmathcalLIn}
\In P_x=\overline\Ad_p\circ\mathcal L_p,\quad p\in P_x,
\end{equation}
does not depend on the choice of $p\in P_x$; here $\overline\Ad_p:\gl(\R^k)/\mathfrak g\to\gl(E_x)/\mathfrak g_x$
is obtained from $\Ad_p:\gl(\R^k)\to\gl(E_x)$ by passing to the quotient. We call $\In P_x$ the {\em inner torsion\/}
of the $G$-structure $P$ at the point $x$ with respect to the connection $\nabla$. Obviously, $\In P=0$ if and only
if $\nabla$ is compatible with $P$. It follows from \eqref{eq:relbaromegaGamma} that if $s:U\to P$ is a smooth
local section with $x\in U$ and $\Gamma$ denotes the Christoffel tensor of $\nabla$ with respect to $s$ then
the inner torsion $\In P_x$ is precisely the composition of the map $\Gamma_x:T_xM\to\gl(E_x)$ with
the quotient map $\gl(E_x)\to\gl(E_x)/\mathfrak g_x$. This observation gives a simple method for computing inner
torsions.

Let us compute inner torsions in some specific examples.
\begin{example}\label{exa:1struc}
If $E$ is trivial and $s:M\to\FR(E)$ is a smooth global frame then $P=s(M)$ is a $G$-structure on $E$ with
$G=\{\Id_{\R^k}\}$. For each $x\in M$, we have $G_x=\{\Id_{E_x}\}$ and
$\mathfrak g_x=\{0\}$; the inner torsion $\In P_x:T_xM\to\gl(E_x)$ is equal to the Christoffel
tensor $\Gamma_x:T_xM\to\gl(E_x)$ of $\nabla$ with respect to $s$.
\end{example}

\begin{example}\label{exa:InsemiRiem}
Let $g$ be a semi-Riemannian structure on $E$ of index $r$, i.e., $g$ is a smooth section of $E^*\otimes E^*$
such that $g_x$ is a nondegenerate symmetric bilinear form on $E_x$ of index $r$, for all $x\in M$;
denote by $\langle\cdot,\cdot\rangle_r$ the {\em standard Minkowski inner product of index $r$\/} in $\R^k$ defined by:
\begin{equation}\label{eq:langlerangler}
\langle v,w\rangle_r=\sum_{i=1}^{k-r}v_iw_i-\hskip-9pt\sum_{i=k-r+1}^k\hskip-9ptv_iw_i,\quad v,w\in\R^k.
\end{equation}
We denote by $\FRo(E)$ the set of all $p\in\FR(E)$ that are linear isometries, so that
$P=\FRo(E)$ is a $G$-structure on $E$, where $G=\Or(k-r,r)$ denotes the Lie group of linear isometries
of $(\R^k,\langle\cdot,\cdot\rangle_r)$.
Let $x\in M$ be fixed. Clearly, $G_x$ is the group of linear isometries of $(E_x,g_x)$
and $\mathfrak g_x$ is the Lie algebra of linear endomorphisms of $E_x$ that are $g_x$-anti-symmetric.
We identify $\gl(E_x)/\mathfrak g_x$ with the space $\Sym(E_x)$ of all linear endomorphisms of $E_x$ that are $g_x$-symmetric
via the map:
\begin{equation}\label{eq:identwithsym}
\gl(E_x)/\mathfrak g_x\ni T+\mathfrak g_x\longmapsto\tfrac12(T+T^*)\in\Sym(E_x),
\end{equation}
where $T^*:E_x\to E_x$ denotes the transpose of $T$ with respect to $g_x$.
Thus, the inner torsion $\In P_x$ is identified with a linear map from $T_xM$ to $\Sym(E_x)$.
Let $s:U\to P$ be a smooth local section with $x\in U$ and let $e,e'\in E_x$ be fixed; consider the local
sections $\epsilon,\epsilon':U\to E$ such that $\epsilon(x)=e$, $\epsilon'(x)=e'$ and such that the $\R^k$-valued maps
$y\mapsto s(y)^{-1}\big(\epsilon(y)\big)$, $y\mapsto s(y)^{-1}\big(\epsilon'(y)\big)$ are constant. Then
$\bfd^s\epsilon=0$, $\bfd^s\epsilon'=0$ and:
\[\nabla_v\epsilon=\Gamma_x(v)\cdot e,\quad
\nabla_v\epsilon'=\Gamma_x(v)\cdot e',\]
for all $v\in T_xM$. Since $s(y):(\R^k,\langle\cdot,\cdot\rangle_r)\to(E_y,g_y)$ is a linear isometry
for all $y\in U$, the real-valued map $g(\epsilon,\epsilon')$ is constant. Thus:
\begin{multline*}
0=v\big(g(\epsilon,\epsilon')\big)=(\nabla_vg)(e,e')+
g_x(\nabla_v\epsilon,e')+g_x(e,\nabla_v\epsilon')\\
=(\nabla_vg)(e,e')+g_x\big(\Gamma_x(v)\cdot e,e'\big)+g_x\big(e,\Gamma_x(v)\cdot e'\big),
\end{multline*}
for all $v\in T_xM$. Then:
\[g_x\big[\big(\Gamma_x(v)+\Gamma_x(v)^*\big)\cdot e,e'\big]=-(\nabla_vg)(e,e')\]
and (using \eqref{eq:identwithsym}):
\[g_x\big(\In P_x(v),\cdot\big)=\tfrac12\,g_x\big[\big(\Gamma_x(v)+\Gamma_x(v)^*\big),\cdot\big]=-\tfrac12\nabla_vg,\]
for all $x\in M$, $v\in T_xM$. Using $g_x$ to identify $\nabla_vg:E_x\times E_x\to\R$ with a linear endomorphism of $E_x$
we obtain:
\[\In P_x(v)=-\tfrac12\nabla_vg.\]
Thus, the inner torsion of $P$ is essentially the covariant derivative of the semi-Riemannian
structure $g$. In particular, $\In P=0$ if and only if $\nabla g=0$.
\end{example}

\begin{example}\label{exa:inntorFREF}
Let $F$ be a vector subbundle of $E$ of rank $l$. For $x\in M$, set:
\[\FR(E_x;F_x)=\big\{p\in\FR(E_x):p(\R^l\oplus\{0\})=F_x\big\}\]
and $\FR(E;F)=\bigcup_{x\in M}\FR(E_x;F_x)$. Then $P=\FR(E;F)$ is a $G$-structure on $E$ with $G$
the Lie subgroup of $\GL(\R^k)$ consisting of linear isomorphisms that preserve $\R^l\oplus\{0\}$.
Let $x\in M$ be fixed. Clearly $G_x$ is the Lie group of linear isomorphisms of $E_x$ that preserve $F_x$
and $\mathfrak g_x$ is the Lie algebra of linear endomorphisms of $E_x$ that preserve $F_x$.
We identify the quotient $\gl(E_x)/\mathfrak g_x$ with the space $\Lin(F_x,E_x/F_x)$ via the map:
\begin{equation}\label{eq:LinEFE}
\gl(E_x)/\mathfrak g_x\ni T+\mathfrak g_x\longmapsto\mathfrak q\circ T\vert_{F_x}\in\Lin(F_x,E_x/F_x),
\end{equation}
where $\mathfrak q:E_x\to E_x/F_x$ denotes the quotient map. Thus, the inner torsion $\In P_x$ is identified
with a linear map from $T_xM$ to $\Lin(F_x,E_x/F_x)$. Let $s:U\to P$ be a smooth local section with $x\in U$
and $e\in F_x$ be fixed. Define $\epsilon$ as in Example~\ref{exa:InsemiRiem}, so that $\nabla_v\epsilon=\Gamma_x(v)\cdot e$, for all $v\in T_xM$.
We have $\epsilon(U)\subset F$ and therefore:
\[\Gamma_x(v)\cdot e+F_x=\nabla_v\epsilon+F_x=\alpha^F_x(v,e)\in E_x/F_x,\]
where $\alpha^F$ denotes the absolute second fundamental form of the vector subbundle $F$.
Hence, using \eqref{eq:LinEFE}:
\[\In P_x(v)=\alpha^F_x(v,\cdot)\in\Lin(F_x,E_x/F_x),\]
for all $x\in M$ and all $v\in T_xM$. In particular, $\In P=0$ if and only if $\alpha^F=0$, i.e., if and only if the covariant
derivative of any smooth section of $F$ is a smooth section of $F$.
\end{example}

\begin{example}\label{exa:calcInPOrEF}
Let $g$ be a semi-Riemannian structure on $E$ of index $r$ and $F$ be a vector subbundle of $E$ such that
the restriction of $g$ to $F$ is a semi-Riemannian structure on $F$ of index $s\le r$; denote
by $F^\perp$ the $g$-orthogonal complement of $F$ in $E$, so that $E=F\oplus F^\perp$. Define $\langle\cdot,\cdot\rangle_r$
as in \eqref{eq:langlerangler} and fix any subspace $F_0$ of $\R^k$ such that the restriction of $\langle\cdot,\cdot\rangle_r$
to $F_0$ is a nondegenerate symmetric bilinear form of index $s$. Then $P=\FRo(E;F)=\bigcup_{x\in M}\FRo(E_x;F_x)$, where
\[\FRo(E_x;F_x)=\big\{p\in\FRo(E_x):p(F_0)=F_x\big\},\]
is a $G$-structure on $E$ where $G$ is the Lie group of linear isometries of $(\R^k,\langle\cdot,\cdot\rangle_r)$
that preserve $F_0$. Denote by $\mathfrak q:E\to F^\perp$ the projection with respect to the decomposition
$E=F\oplus F^\perp$. Let $x\in M$ be fixed.
Clearly $G_x$ is the Lie group of linear isometries of $(E_x,g_x)$ that preserve $F_x$
and $\mathfrak g_x$ is the Lie algebra of $g_x$-anti-symmetric linear endomorphisms of $E_x$ that preserve $F_x$.
We have an isomorphism:
\begin{align*}
\gl(E_x)/\mathfrak g_x&\longrightarrow\Sym(E_x)\oplus\Lin(F_x,F_x^\perp)\\
T+\mathfrak g_x&\longmapsto\big(\tfrac12(T+T^*),\tfrac12\,\mathfrak q_x\circ(T-T^*)\vert_{F_x}\big),
\end{align*}
so that we identify $\In P_x$ with a linear map from $T_xM$ to $\Sym(E_x)\oplus\Lin(F_x,F_x^\perp)$.
Consider the component $\alpha\in\bfG(TM^*\otimes F^*\otimes F^\perp)$ of the connection $\nabla$ with respect to
the decomposition $E=F\oplus F^\perp$. Arguing as in Examples~\ref{exa:InsemiRiem} and \ref{exa:inntorFREF}, one easily computes:
\[\In P_x(v)=\big({-\tfrac12\nabla_vg},\alpha_x(v,\cdot)+\tfrac12\,\mathfrak q\circ\nabla_vg\vert_{F_x}\big),\]
for all $x\in M$, $v\in T_xM$, where $\nabla_vg$ is identified with a linear endomorphism of $E_x$ using $g_x$.
In particular, $\In P=0$ if and only if $\nabla g=0$ and $\alpha=0$, i.e., if and only if
$\nabla g=0$ and the covariant derivative of any smooth section of $F$ is a smooth section of $F$.
\end{example}

\begin{example}\label{exa:FRoe0epsilon}
Let $\epsilon\in\bfG(E)$ be a smooth section of $E$ with $\epsilon(x)\ne0$, for
all $x\in M$. Fix a nonzero vector $e_0\in\R^k$; then:
\[P=\bigcup_{x\in M}\big\{p\in\FR(E_x):p(e_0)=\epsilon(x)\big\}\]
is a $G$-structure on $E$ where $G$ is the subgroup of $\GL(\R^k)$ consisting of isomorphisms that fix $e_0$.
Let $x\in M$ be fixed. Then $G_x$ is the subgroup of $\GL(E_x)$ consisting of isomorphisms that fix $\epsilon(x)$
and $\mathfrak g_x$ is the Lie algebra of linear endomorphisms $T:E_x\to E_x$ such that $T\big(\epsilon(x)\big)=0$.
We identify the quotient $\gl(E_x)/\mathfrak g_x$ with $E_x$ via the map:
\[\gl(E_x)/\mathfrak g_x\ni T+\mathfrak g_x\longmapsto T\big(\epsilon(x)\big)\in E_x,\]
so that $\In P_x$ is identified with a linear map from $T_xM$ to $E_x$.
Let $s:U\to P$ be a smooth local section with $x\in U$. We have $s(y)^{-1}\big(\epsilon(y)\big)=e_0$, for all $y\in U$,
so that $\bfd^s_v\epsilon=0$ and $\nabla_v\epsilon=\Gamma_x(v)\cdot\epsilon(x)$, for all $v\in T_xM$.
Then:
\[\In P_x(v)=\nabla_v\epsilon,\]
for all $v\in T_xM$, i.e., the inner torsion $\In P$ is identified with the covariant derivative of $\epsilon$.
In particular, $\In P=0$ if and only if the section $\epsilon$ is parallel.

Assume now that $g$ is a semi-Riemannian structure on $E$ of index $r$, $\langle\cdot,\cdot\rangle_r$
is defined as in \eqref{eq:langlerangler} and that $g_x\big(\epsilon(x),\epsilon(x)\big)=\langle e_0,e_0\rangle_r$,
for all $x\in M$. Then:
\[P'=\bigcup_{x\in M}\big\{p\in\FRo(E_x):p(e_0)=\epsilon(x)\big\}\]
is a $G$-structure on $E$ where $G$ is the Lie subgroup of $\Or(k-r,r)$ consisting of linear isometries that fix $e_0$.
Let $x\in M$ be fixed. Then $G_x$ is the Lie group of linear isometries of $(E_x,g_x)$ that fix $\epsilon(x)$
and $\mathfrak g_x$ is the Lie algebra of $g_x$-anti-symmetric linear endomorphisms $T$ of $E_x$ such that
$T\big(\epsilon(x)\big)=0$. We have the following linear isomorphism:
\[\gl(E_x)/\mathfrak g_x\ni T+\mathfrak g_x\longmapsto\big(\tfrac12(T+T^*),\tfrac12(T-T^*)\cdot\epsilon(x)\big)
\in\Sym(E_x)\oplus\epsilon(x)^\perp\]
where $\epsilon(x)^\perp$ denotes the $g_x$-orthogonal complement of $\epsilon(x)$ in $E_x$.
Arguing as before, we obtain:
\begin{equation}\label{eq:InPximetric}
\In{P'}_x(v)=\Big({-\tfrac12\nabla_vg},\nabla_v\epsilon+\tfrac12(\nabla_vg)\big(\epsilon(x)\big)\Big),
\end{equation}
for all $x\in M$ and all $v\in T_xM$. In particular, $\In{P'}=0$ if and only if $\nabla g=0$ and $\nabla\epsilon=0$.
\end{example}

\begin{example}\label{exa:computeInPJ}
Assume $k=2l$ and let $J$ be an almost complex structure on $E$, i.e., $J$ is a smooth section of $E^*\otimes E$ and $J_x$
is a complex structure on $E_x$ for all $x\in M$. Consider the complex structure
$J_0:\R^k\cong\R^l\oplus\R^l\ni(v,w)\mapsto(-w,v)$ on $\R^k$ and set $\FRc(E)=\bigcup_{x\in M}\FRc(E_x)$,
where for each $x\in M$, $\FRc(E_x)$ denotes the set of all complex linear isomorphisms $p:(\R^k,J_0)\to(E_x,J_x)$.
Then $P=\FRc(E)$ is a $G$-structure on $E$ where $G=\GL(\R^k,J_0)$ is the Lie group of complex linear isomorphisms
of $(\R^k,J_0)$. Let $x\in M$ be fixed. Then $G_x=\GL(E_x,J_x)$ and $\mathfrak g_x$ is the Lie algebra of complex linear
endomorphisms of $(E_x,J_x)$. We have an isomorphism:
\[\gl(E_x)/\mathfrak g_x\ni T+\mathfrak g_x\longmapsto[T,J_x]\in\overline{\Lin}(E_x,J_x),\]
where $[T,J_x]=T\circ J_x-J_x\circ T$ and $\overline{\Lin}(E_x,J_x)$ denotes the
space of linear maps $T:E_x\to E_x$ such that $T\circ J_x+J_x\circ T=0$.
Let $s:U\to P$ be a smooth local section with $x\in U$ and let $e\in E_x$ be fixed. We define
a local section $\epsilon:U\to E$ as in Example~\ref{exa:InsemiRiem}, so that $\nabla_v\epsilon=\Gamma_x(v)\cdot e$,
for all $v\in T_xM$. Since $U\ni y\mapsto s(y)^{-1}\big(J_y\cdot\epsilon(y)\big)$ is constant,
it follows that $\bfd^s_v\big(J(\epsilon)\big)=0$ and $\nabla_v\big(J(\epsilon)\big)=\Gamma_x(v)\cdot\big(J_x(e)\big)$; then:
\[\Gamma_x(v)\cdot\big(J_x(e)\big)=\nabla_v\big(J(\epsilon)\big)=(\nabla_vJ)(e)+J_x\big(\Gamma_x(v)\cdot e\big),\]
for all $v\in T_xM$. We therefore obtain:
\[\Gamma_x(v)\circ J_x=\nabla_vJ+J_x\circ\Gamma_x(v)\]
and hence:
\[\In P_x(v)=\nabla_vJ,\]
for all $x\in M$ and all $v\in T_xM$. In particular, $\In P=0$ if and only if $J$ is parallel.
\end{example}

\begin{example}\label{exa:InPFRu}
Assume $k=2l$. Let $J$ be an almost complex structure
on $E$, $g$ be a semi-Riemannian structure on $E$ of index $r=2s$, $J_0$ be the complex structure on $\R^k$
considered in Example~\ref{exa:computeInPJ} and $\langle\cdot,\cdot\rangle$ be the nondegenerate symmetric bilinear
form of index $r$ on $\R^k$ defined by:
\[\langle v,w\rangle=\sum_{i=1}^{l-s}v_iw_i-\hskip-9pt\sum_{i=l-s+1}^l\hskip-9ptv_iw_i
+\sum_{i=l+1}^{k-s}v_iw_i-\hskip-9pt\sum_{i=k-s+1}^k\hskip-9ptv_iw_i,\quad v,w\in\R^k;\]
observe that $J_0$ is anti-symmetric with respect to $\langle\cdot,\cdot\rangle$.
Assume that $J_x$ is $g_x$-anti-symmetric for all $x\in M$.
Set $\FRu(E)=\bigcup_{x\in M}\FRu(E_x)$, where for each $x\in M$, $\FRu(E_x)$ is the set of complex linear
isometries from $(\R^k,J_0,\langle\cdot,\cdot\rangle)$ to $(E_x,J_x,g_x)$.
Then $P=\FRu(E)$ is a $G$-structure on $E$ where $G=\Ur(s,l-s)$ is the Lie group of complex linear isometries
of $(\R^k,J_0,\langle\cdot,\cdot\rangle)$. Let $x\in M$ be fixed. Clearly $G_x$ is the Lie group
of complex linear isometries of $(E_x,J_x,g_x)$ and $\mathfrak g_x$ is the Lie algebra of complex linear
$g_x$-anti-symmetric endomorphisms of $(E_x,J_x)$.
We have a linear isomorphism:
\begin{align*}
\gl(E_x)/\mathfrak g_x&\longrightarrow\Sym(E_x)\oplus\overline{\Lin}_\alt(E_x,J_x)\\
T+\mathfrak g_x&\longmapsto\big(\tfrac12(T+T^*),\tfrac12[T-T^*,J_x]\big),
\end{align*}
where $\overline{\Lin}_\alt(E_x,J_x)$ denotes the space of $g_x$-anti-symmetric linear endomorphisms $T$ of $E_x$
such that $T\circ J_x+J_x\circ T=0$. Arguing as in Examples~\ref{exa:InsemiRiem} and \ref{exa:computeInPJ} we obtain:
\[\In P_x(v)=\big({-\tfrac12\nabla_vg},\nabla_vJ-[\nabla_vg,J_x]\big),\]
for all $x\in M$ and all $v\in T_xM$. In particular, $\In P=0$ if and only if both $g$ and $J$ are parallel.
\end{example}

We conclude the section with a technical lemma that will be used later on.
\begin{lem}\label{thm:lemapreCartan}
Let $(M_1,\nabla^1)$, $(M_2,\nabla^2)$ be $n$-dimensional affine manifolds, $G$ be a Lie subgroup of
$\GL(\R^n)$ and $P^1\subset\FR(TM_1)$, $P^2\subset\FR(TM_2)$ be $G$-structures on $M_1$ and $M_2$,
respectively. Assume that for all $x\in M_1$, $y\in M_2$ and for every $G$-structure preserving
map $\sigma:T_xM_1\to T_yM_2$ we have $\In{P^2}_y\circ\sigma=\overline\Ad_\sigma\circ\In{P^1}_x$.
Let $\gamma:I\to M_1$, $\mu:I\to M_2$ be smooth curves and $p:I\to\FR(TM_1)$, $q:I\to\FR(TM_2)$ be horizontal liftings
of $\gamma$ and $\mu$, respectively. For each $t\in I$, set $\sigma(t)=q(t)\circ p(t)^{-1}:T_{\gamma(t)}M_1\to T_{\mu(t)}M_2$.
If $\sigma(t)\big(\gamma'(t)\big)=\mu'(t)$, for all $t\in I$ and if $\sigma(t_0)$ is $G$-structure
preserving for some $t_0\in I$ then $\sigma(t)$ is $G$-structure preserving for all $t\in I$.
\end{lem}
\begin{proof}
By partitioning $I$, we may assume without loss of generality that there are smooth local sections
$s_1:U\to P^1$, $s_2:V\to P^2$ with $\gamma(I)\subset U$, $\mu(I)\subset V$ and such that $U$ is the domain
of a local chart of $M_1$. Let $\omega^i$ denote the connection form of $\FR(TM_i)$ and set $\bar\omega^i=s_i^*\omega$,
$i=1,2$. Since $p$, $q$ are horizontal, \eqref{eq:nabla1p} gives us:
\begin{equation}\label{eq:EDOtildepq}
\tilde p'(t)+\bar\omega^1_{\gamma(t)}\big(\gamma'(t)\big)\circ\tilde p(t)=0,\quad
\tilde q'(t)+\bar\omega^2_{\mu(t)}\big(\mu'(t)\big)\circ\tilde q(t)=0,
\end{equation}
for all $t\in I$, where $\tilde p(t)=s_1\big(\gamma(t)\big)^{-1}\circ p(t)$ and $\tilde q(t)=s_2\big(\mu(t)\big)^{-1}\circ q(t)$.
Now set $L(t)=s_2\big(\mu(t)\big)^{-1}\circ\sigma(t)\circ s_1\big(\gamma(t)\big)=\tilde q(t)\circ\tilde p(t)^{-1}$,
so that $\sigma(t)$ is $G$-structure preserving if and only if $L(t)\in G$. Now \eqref{eq:EDOtildepq} implies:
\begin{equation}\label{eq:EDOL}
L'(t)=L(t)\circ\bar\omega^1_{\gamma(t)}\big(\gamma'(t)\big)-\bar\omega^2_{\mu(t)}\big(\mu'(t)\big)\circ L(t),
\end{equation}
for all $t\in I$. Since $U$ is the domain of a local chart of $M_1$, there exists a smooth time-dependent
vector field $X:I\times U\to TM_1$ in $U$ such that $\gamma$ is an integral curve of $X$ (for instance, let $X(t,x)\in T_xM_1$
be the vector that has the same coordinates as $\gamma'(t)$ in a local chart). We define a smooth time-dependent
vector field $\mathcal X$ in $U\times V\times\GL(\R^n)$ by setting:
\[\mathcal X(t,x,y,A)=\Big(X(t,x),\widehat A\big(X(t,x)\big),A\circ\bar\omega^1_x\big(X(t,x)\big)-
\bar\omega^2_y\big[\widehat A\big(X(t,x)\big)\big]\circ A\Big),\]
for all $t\in I$, $x\in U$, $y\in V$, $A\in\GL(\R^n)$, where $\widehat A:T_xM_1\to T_yM_2$
is defined by $\widehat A=s_2(y)\circ A\circ s_1(x)^{-1}$. Notice that if $A=L(t)$, $x=\gamma(t)$ and $y=\mu(t)$
then $\widehat A=\sigma(t)$, $X(t,x)=\gamma'(t)$ and $\widehat A\big(X(t,x)\big)=\mu'(t)$;
thus, by \eqref{eq:EDOL}, $(\gamma,\mu,L):I\to U\times V\times\GL(\R^n)$ is an integral curve of $\mathcal X$.
We claim that $\mathcal X(t,x,y,A)\in T(M_1\times M_2\times G)$ for all $t\in I$, $x\in M_1$, $y\in M_2$ and all
$A\in G$. Namely, for all $x\in U$, $y\in V$,
$v\in T_xM_1$, $w\in T_yM_2$, we have:
\begin{equation}\label{eq:baromegaiIn}
\bar\omega^1_x(v)+\mathfrak g=\overline\Ad_{s_1(x)}^{\,-1}\big(\In{P_1}_x(v)\big),\quad
\bar\omega^2_y(w)+\mathfrak g=\overline\Ad_{s_2(y)}^{\,-1}\big(\In{P_2}_y(w)\big).
\end{equation}
If $A\in G$ then $\widehat A:T_xM_1\to T_yM_2$ is $G$-structure preserving and therefore:
\begin{equation}\label{eq:relInPi}
\In{P_2}_y\big(\widehat A(v)\big)=\overline\Ad_{\widehat A}\big(\In{P_1}_x(v)\big).
\end{equation}
From \eqref{eq:baromegaiIn} and \eqref{eq:relInPi} it follows that:
\[\Ad_A\big(\bar\omega^1_x(v)\big)-\bar\omega^2_y\big(\widehat A(v)\big)=
A\circ\bar\omega^1_x(v)\circ A^{-1}-\bar\omega^2_y\big(\widehat A(v)\big)\in\mathfrak g.\]
This proves that $A\circ\bar\omega^1_x(v)-\bar\omega^2_y\big(\widehat A(v)\big)\circ A\in T_AG$
and concludes the proof of the claim. We have established that the set of instants $t\in[a,b]$
such that $L(t)\in G$ is open. The conclusion follows from the elementary lemma below.
\end{proof}

\begin{lem}
Let $\mathcal G$ be a Lie group, $G$ be a (not necessarily closed) Lie subgroup of $\mathcal G$ and $\lambda:I\to\mathcal G$
be a continuous curve. If $\lambda^{-1}(G)$ is open in $I$ and nonempty then $\lambda(I)\subset G$.
\end{lem}
\begin{proof}
It suffices to consider the case where $I=[a,b]$ and $\lambda(a)\in G$. Let $t_0\in\left]a,b\right]$ denote the supremum
of the set of all $t\in[a,b]$ such that $\lambda\big([a,t]\big)\subset G$. Clearly $\lambda\big(\left[a,t_0\right[\big)\subset G$
and the conclusion will follow once we show that $\lambda(t_0)\in G$. Let $\varphi:U\subset\mathcal G\to\widetilde U\subset\R^n$
be a local chart with $\lambda(t_0)\in U$ and such that the intersection of any left coset of $G$ with $U$ is a countable
union of {\em slices\/} $U_c=\varphi^{-1}\big(\tilde U\cap\pi^{-1}(c)\big)$, $c\in\R^{n-k}$, where $\pi:\R^n\to\R^{n-k}$
denotes the projection onto the first $n-k$ coordinates, $n=\Dim(\mathcal G)$, $k=\Dim(G)$ (this is a standard
construction for Lie groups; see, for instance, \cite[Chapter 3]{Warner}). For $\varepsilon>0$ sufficiently small,
$\lambda\big([t_0-\varepsilon,t_0]\big)\subset U$ and thus $\lambda\big(\left[t_0-\varepsilon,t_0\right[\big)
\subset U\cap G$ is contained in a countable union of slices; since $(\pi\circ\varphi\circ\lambda)\big(\left[t_0-\varepsilon,t_0\right[\big)$
is a connected countable subset of $\R^{n-k}$, it follows that $\lambda\big(\left[t_0-\varepsilon,t_0\right[\big)$
is contained in a single slice $U_c$. Since $U_c$ is closed in $U$, it follows that $\lambda(t_0)\in U_c\subset G$.
\end{proof}

\end{section}

\begin{section}{Infinitesimally homogeneous affine manifolds with $G$-structure}
\label{sec:infhom}

Let $(M,\nabla)$ be an $n$-dimensional affine manifold, let $G$ be a Lie subgroup of $\GL(\R^n)$ and let
$P\subset\FR(TM)$ be a $G$-structure on $M$. Denote by $\Tor$ and $R$ respectively the torsion and the curvature
tensors of $\nabla$. Recall that, given $x,y\in M$, a map $\sigma:T_xM\to T_yM$
is {\em $G$-structure preserving\/} if $\sigma\circ p\in P_y$ for some (and hence for all) $p\in P_x$.
A smooth map $f:M\to M$ is said to be {\em $G$-structure preserving\/} if $\dd f_x:T_xM\to T_{f(x)}M$ is $G$-structure preserving
for all $x\in M$.

\begin{defin}\label{thm:definfhom}
We say that the triple $(M,\nabla,P)$ is an {\em infinitesimally homogeneous affine manifold with $G$-structure\/}
if for all $x,y\in M$ it is the case that every $G$-structure preserving map $\sigma:T_xM\to T_yM$ relates
$\Tor_x$ with $\Tor_y$, $R_x$ with $R_y$ and $\In P_x$ with $\In P_y$, i.e., $\Tor_y(\sigma\cdot,\sigma\cdot)=\sigma\circ\Tor_x$,
$R_y(\sigma\cdot,\sigma\cdot)=\sigma\circ R(\cdot,\cdot)\circ\sigma^{-1}$ and
$\In P_y\circ\sigma=\overline{\Ad}_\sigma\circ\In P_x$.
\end{defin}
It is easy to see that $(M,\nabla,P)$ is infinitesimally homogeneous if and only if there exist
bilinear maps $\Tor_\oo:\R^n\times\R^n\to\R^n$, $R_\oo:\R^n\times\R^n\to\gl(\R^n)$ and a linear map
$\In{}_\oo:\R^n\to\gl(\R^n)/\mathfrak g$ such that for every $x\in M$, every $p\in P_x$ relates
$\Tor_\oo$ with $\Tor_x$, $R_\oo$ with $R_x$ and $\In{}_\oo$ with $\In P_x$. We will refer to
$\Tor_\oo$, $R_\oo$ and $\In{}_\oo$ collectively as the {\em characteristic tensors\/} of the infinitesimally
homogeneous affine manifold with $G$-structure $(M,\nabla,P)$.

\begin{rem}\label{thm:versions}
Clearly, the characteristic tensors
$\Tor_\oo$, $R_\oo$ and $\In{}_\oo$ are invariant by the action of $G$. This implies that one can induce ``versions''
of the tensors $\Tor_\oo$, $R_\oo$, $\In{}_\oo$ on every vector space endowed with a $G$-structure. More precisely,
let $Z$ be a real $n$-dimensional vector space endowed with a $G$-structure $P_Z$, i.e., an orbit of the right action
of $G$ on $\FR(Z)$. Denote by $G_Z\subset\GL(Z)$ the Lie group of all $G$-structure preserving automorphisms of
$Z$ and by $\mathfrak g_Z\subset\gl(Z)$ its Lie algebra. Given any $p\in P_Z$, there exists a unique triple
of tensors $\Tor_Z:Z\times Z\to Z$, $R_Z:Z\times Z\to\gl(Z)$, $\In{}_Z:Z\to\gl(Z)/\mathfrak g_Z$ that are related
to $\Tor_\oo$, $R_\oo$, $\In{}_\oo$ by $p$. The $G$-invariance of $\Tor_\oo$, $R_\oo$, $\In{}_\oo$ implies
that $\Tor_Z$, $R_Z$, $\In{}_Z$ do not depend on the choice of $p\in P_Z$.
\end{rem}

\begin{defin}
We say that the triple $(M,\nabla,P)$ is {\em locally homogeneous\/} if for every $x,y\in M$ and every $G$-structure preserving
map $\sigma:T_xM\to T_yM$ there exists an open neighborhood $U\subset M$ of $x$, an open neighborhood $V\subset M$
of $y$ and a smooth $G$-structure preserving affine diffeomorphism $f:U\to V$
such that $f(x)=y$ and $\dd f_x=\sigma$. If for every $x,y\in M$ and every $G$-structure preserving map $\sigma:T_xM\to T_yM$
there exists a smooth $G$-structure preserving affine diffeomorphism $f:M\to M$ with $f(x)=y$ and $\dd f_x=\sigma$
we say that the triple $(M,\nabla,P)$ is {\em (globally) homogeneous}.
\end{defin}

Clearly every (locally) homogeneous affine manifold with $G$-structure is also infinitesimally homogeneous. The converse
also holds, as is proven below.
\begin{prop}\label{thm:inflocglobhom}
Let $(M,\nabla,P)$ be an infinitesimally homogeneous affine manifold with $G$-structure. Then $(M,\nabla,P)$
is locally homogeneous. If, in addition, $(M,\nabla)$ is geodesically
complete and $M$ is (connected and) simply-connected then $(M,\nabla,P)$ is globally homogeneous.
\end{prop}
\begin{proof}
Let $x_0,y_0\in M$ be fixed and let $\sigma_0:T_{x_0}M\to T_{y_0}M$ be a $G$-structure preserving map.
We will use the Cartan--Ambrose--Hicks (CAH) theorem (see, for instance, \cite{Wolf}) to produce a
smooth $G$-structure preserving affine diffeomorphism $f:U\to V$ with $U$ an open neighborhood of $x_0$,
$V$ an open neighborhood of $y_0$, $f(x_0)=y_0$ and $\dd f_{x_0}=\sigma_0$. In order to verify the hypotheses
of the CAH theorem, we consider a geodesic $\gamma:[a,b]\to M$ with $\gamma(a)=x_0$
and a geodesic $\mu:[a,b]\to M$ with $\mu(a)=y_0$ and $\mu'(a)=\sigma_0\big(\gamma'(a)\big)$. Let
$p:[a,b]\to\FR(TM)$, $q:[a,b]\to\FR(TM)$ be horizontal liftings of $\gamma$ and $\mu$ respectively such that
$q(a)=\sigma_0\circ p(a)$; set $\sigma(t)=q(t)\circ p(t)^{-1}:T_{\gamma(t)}M\to T_{\mu(t)}M$.
Lemma~\ref{thm:lemapreCartan} implies that $\sigma(t)$ is $G$-structure preserving for all $t\in[a,b]$;
namely, $\sigma(a)=\sigma_0$ is $G$-structure preserving and, since $\gamma'$ and $\mu'$ are parallel,
$\sigma(t)\big(\gamma'(t)\big)=\mu'(t)$, for all $t\in[a,b]$.
Now, by the infinitesimal homogeneity hypothesis, we have that $\sigma(t)$ relates $\Tor_{\gamma(t)}$ with
$\Tor_{\mu(t)}$ and $R_{\gamma(t)}$ with $R_{\mu(t)}$, for all $t\in[a,b]$. The CAH theorem
now gives us a smooth affine diffeomorphism $f:U\to V$ with $f(x_0)=y_0$ and $\dd f_{x_0}=\sigma_0$, where
$U$, $V$ are open neighborhoods of $x_0$ and $y_0$, respectively. The differentials of such map $f$ are given
by the maps $\sigma(t)$ considered above and thus $f$ is $G$-structure preserving. Finally, under the
assumption of simply-connectedness and geodesical completeness of $(M,\nabla)$, one can easily apply the
global version of the CAH theorem to conclude that $M$ is homogeneous.
\end{proof}

\begin{example}\label{exa:constseccurvinfhom}
Let $(M,g)$ be an $n$-dimensional semi-Riemannian manifold with constant sectional curvature equal to $c\in\R$, i.e.:
\[g_x\big(R_x(v,w)v,w\big)=-c\big(g_x(v,v)g_x(w,w)-g_x(v,w)^2\big),\]
for all $x\in M$ and all $v,w\in T_xM$; denote by $r$ the index of $g$. Setting $G=\Or(n-r,r)$ and
$P=\FRo(TM)$ then $(M,\nabla,P)$ is infinitesimally homogeneous, where $\nabla$ denotes the Levi-Civita connection of $g$.
Given an $n$-dimensional real vector space $Z$ endowed with a nondegenerate
symmetric bilinear form $\langle\cdot,\cdot\rangle$ of index $r$ (which amounts to giving a $G$-structure on $Z$) then
the ``versions'' $\Tor_Z$, $R_Z$, $\In{}_Z$ of the characteristic tensors of $(M,\nabla,P)$ are given by
$\Tor_Z=0$, $\In{}_Z=0$ and:
\[R_Z(z_1,z_2)z_3=c\big(\langle z_2,z_3\rangle z_1-\langle z_1,z_3\rangle z_2\big),\]
for all $z_1,z_2,z_3\in Z$.
\end{example}

\begin{example}\label{exa:constholcurvinfhom}
Let $(M,g)$ be an $n$-dimensional semi-Riemannian manifold with Levi-Civita connection $\nabla$,
$J$ be a $\nabla$-parallel $g$-antisymmetric almost complex structure on $M$;
in this case we call $(M,g,J)$ a {\em semi-Riemannian K\"ahler\/} manifold. Assume that $(M,g,J)$ has constant holomorphic
sectional curvature equal to $c\in\R$, i.e., $g_x\big[R_x\big(v,J_x(v)\big)v,J_x(v)\big]=-cg_x(v,v)^2$,
for all $x\in M$ and all $v\in T_xM$. Denoting by $r$ the index of $g$ and setting $G=\Ur\big(\frac12\,r,\frac12(n-r)\big)$,
$P=\FRu(TM)$ (recall Example~\ref{exa:InPFRu}) then $(M,\nabla,P)$ is infinitesimally homogeneous. If $Z$ is a $n$-dimensional
real vector space endowed with a nondegenerate symmetric bilinear form $\langle\cdot,\cdot\rangle$ of index $r$ and a
$\langle\cdot,\cdot\rangle$-antisymmetric complex structure $J_0$ then the ``versions'' of the characteristic tensors of
$(M,\nabla,P)$ are given by $\Tor_Z=0$, $\In{}_Z=0$ and:
\begin{multline*}
R_Z(z_1,z_2)z_3=-\tfrac c4\big(\langle z_1,z_3\rangle z_2-\langle z_2,z_3\rangle z_1
-\langle z_1,J_0(z_3)\rangle J_0(z_2)\\
+\langle z_2,J_0(z_3)\rangle J_0(z_1)-2\langle z_1,J_0(z_2)\rangle J_0(z_3)\big),
\end{multline*}
for all $z_1,z_2,z_3\in Z$.
\end{example}

\begin{example}\label{exa:1strucLie}
Let $\mathcal G$ be an $n$-dimensional Lie group with Lie algebra $\mathbf g$ and let $p_0:\R^n\to\mathbf g$ be a linear
isomorphism. Let $\nabla$ be a left-invariant connection on $\mathcal G$ and set $P=\big\{\dd L_g(1)\circ p_0:g\in\mathcal G\big\}$,
where $L_g:\mathcal G\to\mathcal G$ denotes left multiplication by $g$. Then $P$ is a $G$-structure on $\mathcal G$ with
$G=\{\Id_{\R^n}\!\}$ the trivial group. The affine manifold with $G$-structure $(\mathcal G,\nabla,P)$ is homogeneous, since
the left translations are affine $G$-structure preserving diffeomorphisms. Let $\Gamma:\mathbf g\to\gl(\mathbf g)$ be the linear
map defined by $\Gamma(X)Y=\nabla_XY$, for all $X,Y\in\mathbf g$, where we identify the elements of the Lie algebra $\mathbf g$
with left invariant vector fields in $\mathcal G$. The characteristic tensors of $(\mathcal G,\nabla,P)$ are given by:
\begin{gather*}
\In{}_\oo(v)=p_0^{-1}\circ\Gamma\big(p_0(v)\big)\circ p_0\in\gl(\R^n),\\
\Tor_\oo(v,w)=p_0^{-1}\Big(\Gamma\big(p_0(v)\big)p_0(w)-\Gamma\big(p_0(w)\big)p_0(v)-[p_0(v),p_0(w)]\Big)\in\R^n,\\
R_\oo(v,w)=p_0^{-1}\circ\Big(\big[\Gamma\big(p_0(v)\big),\Gamma\big(p_0(w)\big)\big]-\Gamma\big([p_0(v),p_0(w)]\big)\Big)\circ p_0
\in\gl(\R^n),
\end{gather*}
for all $v,w\in\R^n$. If $Z$ is a real $n$-dimensional vector space then a $G$-structure on $Z$ is just a linear isomorphism
from $\R^n$ to $Z$, which can be used to push-forward $\In{}_\oo$, $\Tor_\oo$, $R_\oo$ respectively to $\In{}_Z$, $\Tor_Z$
and $R_Z$. If $\nabla$ is the Levi-Civita connection of the semi-Riemannian left invariant metric on $\mathcal G$ corresponding
to the nondegenerate symmetric bilinear form $\langle\cdot,\cdot\rangle$ on $\mathbf g$ then $\Gamma$ is given by:
\[\langle\Gamma(X)Y,Z\rangle=\tfrac12\big({-\langle X,[Y,Z]\rangle}+\langle Y,[Z,X]\rangle+\langle Z,[X,Y]\rangle\big),\]
for all $X,Y,Z\in\mathbf g$.
\end{example}

\begin{example}\label{exa:Benoit}
Let $(M,g)$ be an oriented three-dimensional Riemannian manifold with Levi-Civita connection $\nabla$ such that:
\begin{itemize}
\item there exists a Riemannian submersion of $(M,g)$ onto a two-dimensional Riemannian manifold of constant
sectional curvature equal to $\kappa\in\R$;
\item there exists a unitary smooth vector field $\xi$ on $M$, vertical with respect to the submersion, and
a real number $\tau$ such that $\nabla_v\xi=\tau v\times\xi(x)$, for every $x\in M$, $v\in T_xM$, where $\times$
denotes the vector product on $T_xM$ (determined by the given inner product $g_x$ and the given orientation).
\end{itemize}
For instance, all homogeneous three-dimensional Riemannian manifolds having a Lie group of isometries
of dimension greater than or equal to four satisfy the conditions above (see \cite{Benoit}). Let
$P=\bigcup_{x\in M}P_x$, where $P_x$ denotes the set of all positively oriented linear isometries
$p:\R^3\to T_xM$ with $p(1,0,0)=\xi(x)$; then $P$ is a $G$-structure on $M$, where $G\cong\SO(\R^2)$ is
the group of positively oriented linear isometries of $\R^3$ that fix $(1,0,0)$. Using \eqref{eq:InPximetric},
we see that the inner torsion $\In P_x:T_xM\to\Sym(T_xM)\oplus\xi(x)^\perp$ is given by:
\[\In P_x(v)=\big(0,\tau v\times\xi(x)\big),\]
for all $x\in M$, $v\in T_xM$. The curvature tensor $R$ of $\nabla$ can be computed easily and the resulting
formulas show that $R$ can be written only in terms of $\kappa$, $\tau$, $\xi$ and the metric $g$ (see,
for instance, \cite{Benoit}).
Hence $(M,\nabla,P)$ is infinitesimally homogeneous. Given a three dimensional real vector space $Z$ endowed
with a $G$-structure (i.e., an orientation, an inner product $\langle\cdot,\cdot\rangle$ and a unit vector
$\mathfrak z\in Z$) then the ``versions'' of the characteristic tensors of $(M,\nabla,P)$ are given by:
\begin{gather}
\In{}_Z(v)=(0,\tau v\times\mathfrak z)\in\Sym(Z)\oplus\mathfrak z^\perp,\\
\Tor_Z=0,\\
\begin{aligned}
R_Z(v,w)u=(\kappa-&3\tau^2)\big(\langle v,u\rangle w-\langle w,u\rangle v\big)+(\kappa-4\tau^2)
\big(\langle w,\mathfrak z\rangle\langle u,\mathfrak z\rangle v\\
&\langle w,\mathfrak z\rangle\langle v,\mathfrak z\rangle\mathfrak z-\langle v,u\rangle\langle w,\mathfrak z\rangle\mathfrak z
-\langle v,\mathfrak z\rangle\langle u,\mathfrak z\rangle w\big)\in Z,
\end{aligned}
\end{gather}
for all $v,w,u\in Z$.
\end{example}

\begin{example}\label{exa:prodMPs}
Let $(M_1,\nabla_1)$, $(M_2,\nabla_2)$ be affine manifolds, with $\Dim(M_i)=n_i$, $i=1,2$; let $G_i$ be a Lie subgroup
of $\GL(\R^{n_i})$ and let $P_i$ be a $G_i$-structure on $M_i$, $i=1,2$. Denote by $\nabla$ the connection on $M_1\times M_2$
naturally induced by $\nabla_1$, $\nabla_2$; the product $G=G_1\times G_2$ can be identified (diagonally) with a Lie subgroup
of $\GL(\R^{n_1+n_2})$ and the product $P=P_1\times P_2$ is a $G$-structure on $M$. If $(M_1,\nabla_1,P_1)$
and $(M_2,\nabla_2,P_2)$ are homogeneous (resp., infinitesimally homogeneous) then also $(M,\nabla,P)$ is
homogeneous (resp., infinitesimally homogeneous). The characteristic tensors of $(M,\nabla,P)$ are given by the obvious
direct sums of the corresponding characteristic tensors of $(M_1,\nabla_1,P_1)$
and $(M_2,\nabla_2,P_2)$. Notice that if $Z$ is an $(n_1+n_2)$-dimensional real vector space then a $G$-structure on $Z$
is determined by a direct sum decomposition $Z=Z_1\oplus Z_2$ with $\Dim(Z_i)=n_i$, $i=1,2$, by a $G_1$-structure on $Z_1$
and by a $G_2$-structure on $Z_2$.
\end{example}

\end{section}

\begin{section}{Existence of $G$-structure preserving affine immersions}
\label{sec:grotesca}

Throughout the section we consider fixed the following objects: affine manifolds $(M,\nabla)$ and
$(\overline M,\overline\nabla)$, a vector bundle $E^0$ over $M$, a connection $\nabla^0$ on $E^0$
and smooth sections $\alpha^0$, $A^0$ of $TM^*\otimes TM^*\otimes E^0$ and $TM^*\otimes(E^0)^*\otimes TM$ respectively;
we denote by $\Tor$ and $R$ respectively the torsion and the curvature tensors of $\nabla$,
by $\overline\Tor$ and $\overline R$ respectively the torsion and the curvature tensors of $\overline\nabla$,
by $R^0$ the curvature tensor of $\nabla^0$ and by $\nablao$ the connection induced by $\nabla$ and $\nabla^0$
on the vector bundles $TM^*\otimes TM^*\otimes E^0$ and $TM^*\otimes(E^0)^*\otimes TM$.
We are interested in studying conditions under which there is a solution $(f,S)$ for the affine immersion
problem with data $\nabla^0$, $\alpha^0$ and $A^0$, satisfying a given initial condition (i.e., $f(x_0)$,
$\dd f(x_0)$ and $S_{x_0}$ should be prescribed, for some $x_0\in M$).

Set $n=\Dim(M)$, $k=\rank(E^0)$, $\bar n=\Dim(\overline M)$,
\[\widehat E=TM\oplus E^0,\]
and assume that $\bar n=n+k$, i.e., that $\rank(\widehat E)=\bar n$; we denote by $\pi^{TM}:\widehat E\to TM$,
$\pi^{E^0}:\widehat E\to E^0$ the projections. Let $\widehat\nabla$ be the connection in $\widehat E$ whose components are $\nabla$,
$\nabla^0$, $\alpha^0$ and $A^0$. If $(f,S)$ is a solution for the affine immersion
problem with data $\nabla^0$, $\alpha^0$ and $A^0$, we define a vector bundle isomorphism
$L:\widehat E\to f^*T\overline M$ by setting:
\begin{equation}\label{eq:defLfromfS}
L\vert_{TM}=\dd f:TM\to f^*T\overline M,\quad L\vert_{E^0}=S:E^0\to f^*T\overline M.
\end{equation}
The vector bundle isomorphism $L$ is connection preserving. Conversely, given a smooth map $f:M\to\overline M$
and a connection preserving vector bundle isomorphism $L:\widehat E\to f^*T\overline M$ such that
$L\vert_{TM}=\dd f$ then, setting $S=L\vert_{E^0}$, we obtain a solution $(f,S)$ for the affine immersion
problem with data $\nabla^0$, $\alpha^0$ and $A^0$.

Denote by $\omega^M$, $\omega^{\overline M}$, $\theta^M$, $\theta^{\overline M}$ respectively
the $\gl(\R^{\bar n})$-valued connection form of $\widehat\nabla$ in $\FR(\widehat E)$, the
$\gl(\R^{\bar n})$-valued connection form of $\overline\nabla$ in $\FR(T\overline M)$,
the $\R^{\bar n}$-valued $\iota$-canonical form of $\FR(\widehat E)$ and the
$\R^{\bar n}$-valued canonical form of $\FR(T\overline M)$, where $\iota:TM\to\widehat E$ denotes the inclusion.
Let $s:V\to\FR(\widehat E)$ be a smooth local section, where $V$ is an open subset of $M$.
If $f:M\to\overline M$ is a smooth map and $L:\widehat E\to f^*T\overline M$ is a vector bundle isomorphism
then we define a smooth map $F:V\to\FR(T\overline M)$ by setting:
\begin{equation}\label{eq:defFfromfL}
F(x)=L_x\circ s(x),
\end{equation}
for all $x\in V$. It is easy to see that:
\begin{equation}\label{eq:equivalences}
\begin{aligned}
F^*\theta^{\overline M}=s^*\theta^M\quad&\Longleftrightarrow\quad L\vert_{TM}=\dd f,\\
F^*\omega^{\overline M}=s^*\omega^M\quad&\Longleftrightarrow\quad\text{$L$ is connection preserving}.
\end{aligned}
\end{equation}
The conditions $F^*\theta^{\overline M}=s^*\theta^M$ and $F^*\omega^{\overline M}=s^*\omega^M$
can be summarized by writing:
\begin{equation}\label{eq:Fstarlambda}
F^*\lambda^{\overline M}=\lambda^M,
\end{equation}
where:
\begin{equation}\label{eq:deflambdas}
\lambda^{\overline M}=(\theta^{\overline M},\omega^{\overline M}),\quad
\lambda^M=s^*(\theta^M,\omega^M).
\end{equation}
Hence if $(f,S)$ is a solution for the affine immersion problem with data $\nabla^0$, $\alpha^0$ and $A^0$
and if $L$ and $F$ are defined as in \eqref{eq:defLfromfS} and \eqref{eq:defFfromfL} then,
by \eqref{eq:equivalences}, equality \eqref{eq:Fstarlambda} holds.
Notice also that, conversely, if a smooth map $F:V\to f^*T\overline M$ is given then one
can obtain a smooth map $f:V\to\overline M$ and a vector bundle isomorphism
$L:\widehat E\vert_V\to f^*T\overline M$ by setting:
\begin{equation}\label{eq:fLfromF}
f=\overline\pip\circ F,\quad
L_x=F(x)\circ s(x)^{-1}
\end{equation}
for all $x\in V$, where $\overline\pip:\FR(T\overline M)\to\overline M$ denotes the
projection; in particular, if $F$ satisfies \eqref{eq:Fstarlambda}
then one obtain a local solution $(f,S)$ for the affine immersion problem with data $\nabla^0$, $\alpha^0$ and $A^0$
by setting $S=L\vert_{E^0}$.

We prove a general result concerning uniqueness of solutions of the affine immersion problem.
\begin{prop}\label{thm:uniqueaffineimmersion}
Assume that $M$ is connected. If $(f^1,S^1)$, $(f^2,S^2)$ are both solutions of the affine
immersion problem with data $\nabla^0$, $\alpha^0$ and $A^0$ and if there exists $x_0\in M$ with:
\[f^1(x_0)=f^2(x_0),\quad \dd f^1(x_0)=\dd f^2(x_0),\quad S^1_{x_0}=S^2_{x_0},\]
then $(f^1,S^1)=(f^2,S^2)$.
\end{prop}
\begin{proof}
We will show that the set $A$ of points $x\in M$ with $f^1(x)=f^2(x)$, $\dd f^1(x)=\dd f^2(x)$ and
$S^1_x=S^2_x$ is open and closed in $M$ (clearly $x_0\in A$). Let $x\in M$ be fixed and
let $s:V\to\FR(\widehat E)$ be a smooth local section, where $V$ is a connected open neighborhood of $x$ in $M$.
For $i=1,2$, let $L^i:\widehat E\to(f^i)^*T\overline M$, $F^i:V\to\FR(T\overline M)$ be the maps defined from
$(f^i,S^i)$ as in \eqref{eq:defLfromfS} and \eqref{eq:defFfromfL}. As it was observed in \eqref{eq:Fstarlambda}:
\[(F^i)^*\lambda^{\overline M}=\lambda^M.\]
Clearly, for $y\in V$, we have $y\in A$ if and only if $F^1(y)=F^2(y)$. Thus, if
$x$ is not in $A$ then $F^1(x)\ne F^2(x)$ and
therefore $F^1(y)\ne F^2(y)$ for $y$ near $x$. This proves that $A$ is closed.
If $x\in A$ then $F^1(x)=F^2(x)$ and, keeping in mind that $\lambda^{\overline M}_{\bar p}:T_{\bar p}\FR(T\overline M)
\to\R^{\bar n}\oplus\gl(\R^{\bar n})$ is an isomorphism for all $\bar p\in\FR(T\overline M)$, we apply Lemma~\ref{thm:uniquesolutionPDE}
below and obtain that $F^1=F^2$. Thus $V\subset A$ and we are done.
\end{proof}

\begin{lem}\label{thm:uniquesolutionPDE}
Let $\mathcal M$, $\mathcal N$ be differentiable manifolds, $Z$ be a real finite-dimensional vector space and
$\lambda$, $\tilde\lambda$, be $Z$-valued smooth $1$-forms on $\mathcal M$ and on $\mathcal N$
respectively; assume that $\mathcal M$ is connected and that $\tilde\lambda_y:T_y\mathcal N\to Z$ is an
isomorphism, for all $y\in\mathcal N$. If $F_1:\mathcal M\to\mathcal N$,
$F_2:\mathcal M\to\mathcal N$ are smooth maps with $F_i^*\tilde\lambda=\lambda$, $i=1,2$, and if
$F_1(x_0)=F_2(x_0)$ for some $x_0\in\mathcal M$ then $F_1=F_2$.
\end{lem}
\begin{proof}
If $\gamma:I\to\mathcal M$ is a smooth curve such that $F_1\big(\gamma(t_0)\big)=F_2\big(\gamma(t_0)\big)$ for some
$t_0\in I$ then $F_1\circ\gamma=F_2\circ\gamma$; namely, both $F_1\circ\gamma$ and $F_2\circ\gamma$ are integral
curves of the smooth time-dependent vector field:
\[I\times\mathcal N\ni(t,y)\longmapsto\tilde\lambda_y^{-1}\big[\lambda_{\gamma(t)}\big(\gamma'(t)\big)\big]
\in T_y\mathcal N.\]
The conclusion follows from the observation that, since $\mathcal M$ is connected, any two points of $\mathcal M$ can be joined
by a piecewise smooth curve.
\end{proof}

In the remainder of this section we will also fix a Lie subgroup $G$ of $\GL(\R^{\bar n})$
with Lie algebra $\mathfrak g\subset\gl(\R^n)$, a $G$-structure $\widehat P$ on $\widehat E$ and
a $G$-structure $\overline P$ on $\overline M$.

\begin{defin}\label{thm:LdftogetherS}
A (possibly local) solution $(f,S)$ for the affine
immersion problem with data $\nabla^0$, $\alpha^0$ and $A^0$ is said to be {\em $G$-structure preserving\/}
if the vector bundle isomorphism $L:\widehat E\to f^*T\overline M$ defined in \eqref{eq:defLfromfS}
is $G$-structure preserving, i.e., if $L_x\circ p\in\overline P_{f(x)}$, for all $x\in M$ and all
$p\in\widehat P_x$.
\end{defin}

If $(\overline M,\overline\nabla,\overline P)$ is infinitesimally homogeneous
then every $\bar n$-dimensional real vector space $Z$ endowed with a $G$-structure
inherits ``versions'' $\overline\Tor_Z:Z\times Z\to Z$, $\overline R_Z:Z\times Z\to\gl(Z)$
and $\overline{\In{}}_Z:Z\to\gl(Z)/\mathfrak g_Z$ of the characteristic tensors $\overline\Tor_\oo:\R^{\bar n}\times\R^{\bar n}\to\R^{\bar n}$,
$\overline R_\oo:\R^{\bar n}\times\R^{\bar n}\to\gl(\R^{\bar n})$ and
$\overline{\In{}}_\oo:\R^{\bar n}\to\gl(\R^{\bar n})/\mathfrak g$ of $(\overline M,\overline\nabla,\overline P)$
(recall Remark~\ref{thm:versions}).
\begin{teo}\label{thm:grotescao}
Assume that $(\overline M,\overline\nabla,\overline P)$
is infinitesimally homogeneous with characteristic tensors $\overline\Tor_\oo$,
$\overline R_\oo$ and $\overline{\In{}}_\oo$ and that the equalities:
\begin{gather}
\label{eq:Gausschar}
\begin{aligned}
\pi^{TM}\big(\overline R_{\widehat E_x}(v,w)u\big)=R_x(v,w)u&+A^0_x\big(v,\alpha^0_x(w,u)\big)\\
&-A^0_x\big(w,\alpha^0_x(v,u)\big),
\end{aligned}\\[5pt]
\label{eq:Riccichar}
\begin{aligned}
\pi^{E^0}\big(\overline R_{\widehat E_x}(v,w)e\big)=R^0_x(v,w)&+\alpha^0_x\big(v,A^0_x(w,e)\big)\\
&-\alpha^0_x\big(w,A^0_x(v,e)\big),
\end{aligned}\\[5pt]
\label{eq:Codazzichar}
\begin{aligned}
\pi^{E^0}\big(\overline R_{\widehat E_x}(v,w)u\big)=(\nablao\alpha^0)_x(v,w,u)&-(\nablao\alpha^0)_x(w,v,u)\\
&+\alpha^0_x\big(\Tor_x(v,w),u\big),
\end{aligned}\\[5pt]
\label{eq:Codazzi2char}
\begin{aligned}
\pi^{TM}\big(\overline R_{\widehat E_x}(v,w)e\big)=(\nablao A^0)_x(v,w,e)
&-(\nablao A^0)_x(w,v,e)\\
&+A^0_x\big(\Tor_x(v,w),e\big),
\end{aligned}\\[5pt]
\label{eq:torsioneq1char}\pi^{TM}\big(\overline\Tor_{\widehat E_x}(v,w)\big)=\Tor_x(v,w),\\
\label{eq:torsioneq2char}\pi^{E^0}\big(\overline\Tor_{\widehat E_x}(v,w)\big)=\alpha^0_x(v,w)-\alpha^0_x(w,v),\\
\label{eq:inntorchar}\overline{\In{}}_{\widehat E_x}(v)=\In{\widehat P}_x(v),
\end{gather}
hold, for all $x\in M$, $v,w,u\in T_xM$ and all $e\in E^0_x$. Then, for every $x_0\in M$, every $y_0\in\overline M$
and every $G$-structure preserving linear map $\sigma_0:\widehat E_{x_0}\to T_{y_0}\overline M$ there
exists a $G$-structure preserving local solution $(f,S)$ for the affine
immersion problem with data $\nabla^0$, $\alpha^0$ and $A^0$ whose domain is an open neighborhood
$U$ of $x_0$ in $M$, $f(x_0)=y_0$ and $\sigma_0=L_{x_0}$, where $L$ is as in \eqref{eq:defLfromfS}.
Moreover, if $M$ is (connected and) simply-connected and $(\overline M,\overline\nabla)$ is geodesically complete
then there exists a unique global solution $(f,S)$ for the affine immersion problem with data
$\nabla^0$, $\alpha^0$ and $A^0$ such that $f(x_0)=y_0$ and $\sigma_0=L_{x_0}$.
\end{teo}
\begin{proof}[Proof of the existence of local solutions]
Let $s:V\to\widehat P$ be a smooth local section defined in an open neighborhood $V$ of $x_0$ in $M$
and consider the $\gl(\R^{\bar n})$-valued $1$-forms $\lambda^M$ and $\lambda^{\overline M}$ defined
in \eqref{eq:deflambdas}. We will look for a smooth map $F:U\to\widehat P\subset\FR(T\overline M)$
such that $F^*\lambda^{\overline M}=\lambda^M\vert_U$ and $F(x_0)=\sigma_0\circ s(x_0)$, where $U$ is an open neighborhood of $x_0$
in $V$. Once such map $F$ is found, we define $f:U\to\overline M$ and $L:\widehat E\vert_U\to f^*T\overline M$
as in \eqref{eq:fLfromF} and we set $S=L\vert_{E^0}$. It will then follow from \eqref{eq:equivalences}
that $(f,S)$ is a local solution for the affine immersion problem with data
$\nabla^0$, $\alpha^0$ and $A^0$; moreover, the fact that $s$ takes values in $\widehat P$ and
$F$ takes values in $\overline P$ will imply that $L$ is $G$-structure preserving.

In order to find the map $F$ we will use the version of the Frobenius theorem stated in Lemma~\ref{thm:formsFrobenius}
below. We claim that for each $y\in\overline M$ and each $\bar p\in\overline P_y$,
the linear map $\lambda^{\overline M}_{\bar p}:T_{\bar p}\overline P\to\R^{\bar n}\oplus\gl(\R^{\bar n})$ is an isomorphism onto the space:
\[Z=\big\{(u,X)\in\R^{\bar n}\oplus\gl(\R^{\bar n}):\overline{\In{}}_\oo(u)=X+\mathfrak g\big\}.\]
Namely, $(\theta^{\overline M}_{\bar p},\omega^{\overline M}_{\bar p})=(\bar p^{-1}\oplus\Id)\circ
(\dd\overline\pip_{\bar p},\omega^{\overline M}_{\bar p})$ and the isomorphism $(\dd\overline\pip_{\bar p},\omega^{\overline M}_{\bar p})$
carries $T_{\bar p}\overline P$ onto the space (recall \eqref{eq:defcalVp} and \eqref{eq:relmathcalLIn}):
\[\mathcal V_{\bar p}=\big\{(v,X)\in T_y\overline M\oplus\gl(\R^{\bar n}):
\overline\Ad_{\bar p}^{\,-1}\big(\In{\overline P}_y(v)\big)=X+\mathfrak g\big\},\]
and $\overline{\In{}}_\oo=\overline\Ad_{\bar p}^{\,-1}\circ\In{\overline P}_y\circ\bar p$.
We now claim that for all $x\in V$, the linear map $\lambda^M_x:T_xM\to\R^{\bar n}\oplus\gl(\R^{\bar n})$ takes
values in $Z$. Namely, for all $\hat p\in\widehat P_x$,
the linear map $(\theta^M_{\hat p},\omega^M_{\hat p})=\big((\hat p^{\,-1}\circ\iota_x)\oplus\Id\big)\circ
(\dd\widehat\pip_{\hat p},\omega^M_{\hat p}):T_{\hat p}\FR(\widehat E)\to\R^{\bar n}\oplus\gl(\R^{\bar n})$
takes $T_{\hat p}\widehat P$ to the space:
\begin{equation}\label{eq:espacocomplicado}
\big\{(u,X)\in\R^{\bar n}\oplus\gl(\R^{\bar n}):
\hat p(u)\in T_xM,\ \overline\Ad_{\hat p}^{\,-1}\big[\In{\widehat P}_x\big(\hat p(u)\big)\big]=X+\mathfrak g\big\},
\end{equation}
where $\widehat\pip:\FR(\widehat E)\to M$ denotes the projection. By assumption \eqref{eq:inntorchar}
we may replace $\In{\widehat P}_x$ with $\overline{\In{}}_{\widehat E_x}$ in \eqref{eq:espacocomplicado} and,
since $\overline\Ad_{\hat p}^{\,-1}\circ\overline{\In{}}_{\widehat E_x}\circ\hat p=\overline{\In{}}_\oo$,
we conclude that $(\theta^M_{\hat p},\omega^M_{\hat p})(T_{\hat p}\widehat P)\subset Z$. The proof of the
claim is completed by setting $\hat p=s(x)$ and by observing that $\lambda^M_x(T_xM)$ is contained
in $(\theta^M_{\hat p},\omega^M_{\hat p})(T_{\hat p}\widehat P)$. We have so far checked the validity of the
hypotheses of Lemma~\ref{thm:formsFrobenius}; we will now verify that condition (b) in the statement of
the lemma holds and this will give us the desired map $F:U\to\overline P$ with $F^*\lambda^{\overline M}=\lambda^M\vert_U$
and $F(x_0)=\sigma_0\circ s(x_0)$. Let $x\in V$, $y\in\overline M$, $\bar p\in\overline P_y$ be fixed and
set $\tau=(\lambda^{\overline M}_{\bar p})^{-1}\circ\lambda^M_x:T_xM\to T_{\bar p}\overline P$. We have to check
that:
\begin{equation}\label{eq:wehavetocheck}
\tau^*\dd\theta^{\overline M}_{\bar p}=\dd(s^*\theta^M)_x,\quad
\tau^*\dd\omega^{\overline M}_{\bar p}=\dd(s^*\omega^M)_x;
\end{equation}
since $\tau^*\lambda^{\overline M}_{\bar p}=\lambda^M_x$, i.e., $\tau^*\theta^{\overline M}_{\bar p}=s^*\theta^M_x$
and $\tau^*\omega^{\overline M}_{\bar p}=s^*\omega^M_x$, equalities \eqref{eq:defThetaOmega} imply that
\eqref{eq:wehavetocheck}
is equivalent to:
\begin{equation}\label{eq:checkisequivalentto}
\tau^*\Theta^{\overline M}_{\bar p}=(s^*\Theta^M)_x,\quad
\tau^*\Omega^{\overline M}_{\bar p}=(s^*\Omega^M)_x,
\end{equation}
where $\Omega^M$, $\Omega^{\overline M}$, $\Theta^M$, $\Theta^{\overline M}$ denote respectively
the curvature form of $\widehat\nabla$ in $\FR(\widehat E)$, the curvature form of $\overline\nabla$
in $\FR(T\overline M)$, the $\iota$-torsion form of $\widehat\nabla$ in $\FR(\widehat E)$ and the
curvature form of $\overline\nabla$ in $\FR(T\overline M)$. It is easy to see that
$\dd\overline\pip_{\bar p}\circ\tau=\bar p\circ s(x)^{-1}\vert_{T_xM}$; using this equality,
and equalities \eqref{eq:ThetaT} and \eqref{eq:OmegaR}, we get that \eqref{eq:checkisequivalentto} is equivalent to:
\begin{gather}
\label{eq:equivT1}\bar p^{\,-1}\Big(\overline\Tor_y\big[\bar p\big(s(x)^{-1}(v),s(x)^{-1}(w)\big)\big]\Big)
=s(x)^{-1}\big(\widehat\Tor^\iota_x(v,w)\big),\\
\label{eq:equivR1}\bar p^{\,-1}\circ\overline R_y\big[\bar p\big(s(x)^{-1}(v),s(x)^{-1}(w)\big)\big]\circ\bar p
=s(x)^{-1}\circ\widehat R_x(v,w)\circ s(x),
\end{gather}
for all $v,w\in T_xM$, where $\widehat\Tor^\iota$ and $\widehat R$ denote respectively the $\iota$-torsion tensor and
the curvature tensor of $\widehat\nabla$. Clearly \eqref{eq:equivT1} and \eqref{eq:equivR1} are equivalent to:
\begin{gather}
\label{eq:equivT2}\overline\Tor_\oo\big(s(x)^{-1}(v),s(x)^{-1}(w)\big)=s(x)^{-1}\big(\widehat\Tor^\iota_x(v,w)\big),\\
\label{eq:equivR2}\overline R_\oo\big(s(x)^{-1}(v),s(x)^{-1}(w)\big)
=s(x)^{-1}\circ\widehat R_x(v,w)\circ s(x).
\end{gather}
Now \eqref{eq:equivT2} and \eqref{eq:equivR2} are equivalent to:
\begin{gather}
\label{eq:equivT3}\overline\Tor_{\widehat E_x}(v,w)=\widehat\Tor^\iota_x(v,w),\quad v,w\in T_xM,\\
\label{eq:equivR3}\overline R_{\widehat E_x}(v,w)=\widehat R_x(v,w),\quad v,w\in T_xM.
\end{gather}
Finally, \eqref{eq:equivT3} is equivalent to assumptions \eqref{eq:torsioneq1char}, \eqref{eq:torsioneq2char}
(recall \eqref{eq:gentorsioneq1} and \eqref{eq:gentorsioneq2}), while \eqref{eq:equivR3} is equivalent to assumptions
\eqref{eq:Gausschar}, \eqref{eq:Riccichar}, \eqref{eq:Codazzichar} and \eqref{eq:Codazzi2char}
(recall \eqref{eq:genGauss}, \eqref{eq:genRicci}, \eqref{eq:genCodazzi} and \eqref{eq:genCodazzi2}). This concludes the
proof of the existence of local solutions of the affine immersion problem.
\end{proof}

\begin{lem}\label{thm:formsFrobenius}
Let $\mathcal M$, $\mathcal N$ be differentiable manifolds, $Z$ be a real finite-dimensional vector space and
$\lambda$, $\tilde\lambda$, be $Z$-valued smooth $1$-forms on $\mathcal M$ and on $\mathcal N$
respectively; assume that $\tilde\lambda_y:T_y\mathcal N\to Z$ is an isomorphism, for all $y\in\mathcal N$.
The following conditions are equivalent:
\begin{itemize}
\item[(a)] for all $x\in\mathcal M$, $y\in\mathcal N$ there exists a smooth map $F:U\to\mathcal N$ defined in an
open neighborhood $U$ of $x$ in $\mathcal M$ with $F(x)=y$ such that $F^*\tilde\lambda=\lambda\vert_U$;
\item[(b)] for all $x\in\mathcal  M$, $y\in\mathcal N$, $\tau_{xy}^*\dd\tilde\lambda_y=\dd\lambda_x$,
where $\tau_{xy}:T_x\mathcal M\to T_y\mathcal N$ is the linear map defined by
$\tau_{xy}=\tilde\lambda_y^{-1}\circ\lambda_x$.
\end{itemize}
\end{lem}
\begin{proof}
It is a simple application of the classical Frobenius theorem to the distribution on $\mathcal M\times\mathcal  N$ annihilated by
the $Z$-valued $1$-form $\theta=\pr_2^*\tilde\lambda-\pr_1^*\lambda$, where
$\pr_1:\mathcal M\times\mathcal N\to\mathcal M$, $\pr_2:\mathcal M\times\mathcal N\to\mathcal N$
denote the projections.
\end{proof}

In order to prove the existence of global solutions to the affine immersion problem, we will employ a very general
globalization technique that is explained below. Let $X$ be a topological space and let $\mathfrak P$ be a pre-sheaf
of sets over $X$, i.e., a cofunctor from the category of open subsets of $X$ (partially ordered by inclusion) to the category
of sets and maps; for each open subset $U$ of $X$ we denote by $\mathfrak P(U)$ the corresponding set and
for every open subset $V$ of $U$ we denote by $\mathfrak P_{U,V}:\mathfrak P(U)\to\mathfrak P(V)$ the corresponding
map. We say that the pre-sheaf $\mathfrak P$ has the {\em localization property\/} if,
given a family $(U_i)_{i\in I}$ of open subsets of $X$ and setting $U=\bigcup_{i\in I}U_i$
then the map:
\begin{equation}\label{eq:PUUi}
\mathfrak P(U)\ni\mathfrak f\longmapsto\big(\mathfrak P_{U,U_i}(\mathfrak f)\big)_{i\in I}\in\prod_{i\in I}
\mathfrak P(U_i)
\end{equation}
is injective and its image consists of the families $(\mathfrak f_i)_{i\in I}$ in $\prod_{i\in I}\mathfrak P(U_i)$
such that $\mathfrak P_{U_i,U_i\cap U_j}(\mathfrak f_i)=\mathfrak P_{U_j,U_i\cap U_j}(\mathfrak f_j)$, for all
$i,j\in I$. We say that the pre-sheaf $\mathfrak P$ has the {\em uniqueness property\/}
if for every connected open subset $U\subset X$ and every nonempty open subset $V\subset U$
the map $\mathfrak P_{U,V}$ is injective. We say that an open subset $U\subset X$ has the {\em extension property with respect to the
pre-sheaf $\mathfrak P$\/} if for every connected nonempty open subset $V$ of $U$ the map $\mathfrak P_{U,V}$
is surjective. We say that the pre-sheaf $\mathfrak P$ has the {\em extension property\/}
if $X$ can be covered by open sets having the extension property with respect to $\mathfrak P$.

If $\widetilde X$, $X$ are topological spaces and $\pi:\widetilde X\to X$ is a local homeomorphism then an open
subset $U$ of $X$ is said to be {\em fundamental\/} for $\pi$ if $\pi^{-1}(U)$ is a disjoint union $\bigcup_{i\in I}U_i$
of open subsets $U_i$ of $\widetilde X$ such that $\pi$ maps $U_i$ homeomorphically onto $U$ for all $i\in I$.
An open subset $U$ of $X$ is said to be {\em quasi-fundamental\/} for $\pi$ if for every $x\in U$
and every $\tilde x\in\pi^{-1}(x)$ there exists a continuous local section $s:U\to\widetilde X$ of $\pi$ such that
$s(x)=\tilde x$. Clearly, every fundamental open subset is also quasi-fundamental; conversely, if
$\widetilde X$ is Hausdorff then it is easy to check that every connected quasi-fundamental open subset
of $X$ is fundamental (take the $U_i$'s to be the images of the continuous local sections of $\pi$ defined in $U$).

We have the following general globalization principle:
\begin{prop}\label{thm:bigextensionprop}
Let $X$ be a topological space and $\mathfrak P$ be a pre-sheaf of sets over $X$.
Assume that $X$ is Hausdorff, locally arc-connected, connected, simply-connected and that the pre-sheaf $\mathfrak P$
has the localization property, the uniqueness property and the extension property. Then for every connected nonempty
open subset $V\subset X$ and every $\mathfrak f\in\mathfrak P(V)$ there exists a unique
$\bar{\mathfrak f}\in\mathfrak P(X)$ with $\mathfrak P_{X,V}(\bar{\mathfrak f})=\mathfrak f$.
\end{prop}
\begin{proof}
Let $\pi:\mathcal S\to X$ be the sheaf of germs associated to the pre-sheaf $\mathfrak P$ and for each open subset $U$
of $X$, each $\mathfrak f\in\mathfrak P(U)$ and each $x\in U$, denote by $[\mathfrak f]_x\in\mathcal S_x=\pi^{-1}(x)$
the germ of $\mathfrak f$ at $x$. The fact that $\mathfrak P$
has the localization property means that for every open subset $U$ of $X$ and every
continuous local section $s:U\to\mathcal S$ of $\pi$ there exists a unique
$\mathfrak f\in\mathfrak P(U)$ such that $s(x)=[\mathfrak f]_x$, for all $x\in U$. Since $X$ is locally connected,
it follows that every open subset of $X$ having the extension property with respect to $\mathfrak P$ is quasi-fundamental
for $\pi$. Moreover, since $X$ is Hausdorff and locally connected, the uniqueness property of $\mathfrak P$ implies
that $\mathcal S$ is Hausdorff. Thus, every quasi-fundamental connected open subset of $X$ is fundamental and we obtain
that $\pi$ is a covering map. Let $x\in V$ be fixed and let $\mathcal S'$ denote the connected component of $\mathcal S$
containing $[\mathfrak f]_x$; since $X$ is locally arc-connected, connected and simply-connected, it follows
that the restriction of $\pi$ to $\mathcal S'$ is a homeomorphism onto $X$. The inverse of $\pi\vert_{\mathcal S'}$
is a continuous global section $s:X\to\mathcal S$ of $\pi$ with $s(x)=[\mathfrak f]_x$. As observed in the beginning
of the proof, there exists a unique $\bar{\mathfrak f}\in\mathfrak P(X)$ such that $[\bar{\mathfrak f}]_y=s(y)$,
for all $y\in X$; in particular, we have $[\bar{\mathfrak f}]_x=[\mathfrak f]_x$. Since $V$ is connected and
$\mathfrak P$ has the uniqueness property, it follows that $\mathfrak P_{X,V}(\bar{\mathfrak f})=\mathfrak f$.
\end{proof}

We can now prove the global part of the statement of Theorem~\ref{thm:grotescao}.
\begin{proof}[Proof of the existence of global solutions]
By considering the universal covering, there is no loss of generality in assuming that $\overline M$ is connected and
simply-connected; thus, it follows from Proposition~\ref{thm:inflocglobhom} that $(\overline M,\overline\nabla,\overline P)$
is homogeneous.

For each open subset $U$ of $M$ let $\mathfrak P(U)$ be the set of all $G$-structure preserving local solutions for the affine
immersion problem with data $\nabla^0$, $\alpha^0$, $A^0$ and with domain $U$; if $V$ is an open subset of $U$,
we let $\mathfrak P_{U,V}:\mathfrak P(U)\to\mathfrak P(V)$ be the restriction map. It is obvious that the
pre-sheaf of sets $\mathfrak P$ over $M$ has the localization property and it follows from Proposition~\ref{thm:uniqueaffineimmersion}
that $\mathfrak P$ has the uniqueness property. We will prove that $\mathfrak P$ has the extension property and
then the aimed result will follow from the local part of the statement
of Theorem~\ref{thm:grotescao} and from Proposition~\ref{thm:bigextensionprop}. In order to prove that $\mathfrak P$
has the extension property we will show that every open subset $U$ of $M$ such that $\mathfrak P(U)$ is nonempty
has the extension property for $\mathfrak P$ (observe also that the local part of the statement of Theorem~\ref{thm:grotescao}
implies that every point of $M$ has an open neighborhood $U$ such that $\mathfrak P(U)$ is nonempty).
Let $U$ be an open subset of $M$ such that $\mathfrak P(U)$ is nonempty, $V$ be a connected nonempty open subset
of $U$ and let $(f,S)$ be an element of $\mathfrak P(V)$. We have to find $(\bar f,\overline S)$ in $\mathfrak P(U)$
such that $(f,S)$ is the restriction of $(\bar f,\overline S)$ to $V$. Let $(\tilde f,\widetilde S)$ be an arbitrary
element of the nonempty set $\mathfrak P(U)$ and let $x\in V$ be an arbitrary point. Denote by
$L:\widehat E\vert_V\to f^*T\overline M$, $\widetilde L:\widehat E\vert_U\to\tilde f^*T\overline M$ the
$G$-structure preserving connection preserving vector bundle isomorphisms corresponding respectively
to $(f,S)$ and $(\tilde f,\widetilde S)$, as in \eqref{eq:defLfromfS}. Let $\tau:T_{\tilde f(x)}\overline M\to T_{f(x)}\overline M$
be the linear isomorphism such that $\tau\circ\widetilde L_x=L_x$; clearly, $\tau$ is $G$-structure preserving and
thus there exists a global $G$-structure preserving affine diffeomorphism $\phi:\overline M\to\overline M$ such that
$\phi\big(\tilde f(x)\big)=f(x)$ and $\dd\phi_{\tilde f(x)}=\tau$. We set $\bar f=\phi\circ\tilde f:U\to\overline M$
and we consider the $G$-structure preserving connection preserving vector bundle isomorphism $\overline L:\widehat E\vert_U\to\bar f^*T\overline M$
defined by $\overline L_y=\dd\phi_{\tilde f(y)}\circ\widetilde L_y$, for all $y\in U$. Setting
$\overline S=\overline L\vert_{E^0}$ then $(\bar f,\overline S)$ is in $\mathfrak P(U)$.
Moreover, since $\bar f(x)=f(x)$ and $\overline L_x=L_x$, it follows from Proposition~\ref{thm:uniqueaffineimmersion}
and from the connectedness of $V$ that $(f,S)$ is the restriction of $(\bar f,\overline S)$ to $V$.
This concludes the proof.
\end{proof}

\begin{example}\label{exa:affimmgroup}
Assume that $\overline M$ is a Lie group $\mathcal G$ with Lie algebra $\mathbf g$
and that $\overline\nabla$ is a left invariant connection corresponding to a linear map
$\Gamma:\mathbf g\to\gl(\mathbf g)$ as in Example~\ref{exa:1strucLie}; let $p_0:\R^{\bar n}\to\mathbf g$
be a linear isomorphism and let $\overline P$ be the $\{\Id_{\R^n}\!\}$-structure on $\mathcal G$ defined
by $\overline P=\big\{\dd L_g(1)\circ p_0:g\in\mathcal G\big\}$. Then $(\overline M,\overline\nabla,\overline P)$
is homogeneous. Assume that we are given a global smooth frame $s:M\to\FR(\widehat E)$ of $\widehat E$,
so that $\widehat P=s(M)$ is a $\{\Id_{\R^n}\!\}$-structure on $\widehat E$. An application of
Theorem~\ref{thm:grotescao} shows that, under assumptions \eqref{eq:Gausschar}, \eqref{eq:Riccichar}, \eqref{eq:Codazzichar},
\eqref{eq:Codazzi2char}, \eqref{eq:torsioneq1char}, \eqref{eq:torsioneq2char}, \eqref{eq:inntorchar}, we can
find local solutions $(f,S)$ (with adequately prescribed initial conditions $f(x_0)$, $\dd f(x_0)$, $S_{x_0}$) for the affine
immersion problem with data $\nabla^0$, $\alpha^0$ and $A^0$, where the vector bundle isomorphism
$L:\widehat E\to f^*T\mathcal G$ defined by \eqref{eq:defLfromfS} satisfies $L_x\circ s(x)=\dd L_{f(x)}(1)\circ p_0$,
for all $x\in M$. Observe that the tensors $\overline R_{\widehat E_x}$, $\overline\Tor_{\widehat E_x}$,
$\overline{\In{}}_{\widehat E_x}$ that appear in the lefthand side of the assumptions are the push-forward by
the linear isomorphism $s(x):\R^{\bar n}\to\widehat E_x$ of the respective characteristic tensors $R_\oo$,
$\Tor_\oo$ and $\In{}_\oo$ that appear in Example~\ref{exa:1strucLie}. Equation \eqref{eq:inntorchar} means
simply that the Christoffel tensor $\widehat\Gamma_x:T_xM\to\gl(\widehat E_x)$ of $\widehat\nabla$ with
respect to the frame $s$ is the restriction to $T_xM$ of the push-forward to $\widehat E_x$
of $\In{}_\oo$ by $s(x)$.
\end{example}

\end{section}

\begin{section}{Existence of $G$-structure preserving isometric immersions}

As a special case of the theory developed in Section~\ref{sec:grotesca}, in this section we will consider semi-Riemannian
manifolds endowed with $G$-structures contained in the orthonormal frame bundle, aiming at an existence theorem
for isometric immersions. More precisely, we consider fixed the following objects:
an $\bar n$-dimensional semi-Riemann\-ian manifold
$(\overline M,\bar g)$, where the semi-Riemannian metric $\bar g$ has index $\bar r$,
an $n$-dimensional semi-Riemannian manifold $(M,g)$, where the semi-Riemann\-ian metric $g$
has index $r$, a vector bundle $E^0$ of rank $k$ over $M$ endowed with a semi-Riemannian structure $g^0$ of index
$s$, where $\bar n=n+k$ and $\bar r=r+s$, a connection $\nabla^{E^0}$ on $E^0$ compatible with $g^0$,
a smooth symmetric section $\alpha^0$ of $TM^*\otimes TM^*\otimes E^0$, a Lie subgroup $G$ of $\Or(\bar n-\bar r,\bar r)$
with Lie algebra $\mathfrak g\subset\so(\bar n-\bar r,\bar r)$,
a $G$-structure $\overline P$ on $\overline M$ contained in the orthonormal frame bundle $\FRo(T\overline M)$ and
a $G$-structure $\widehat P$ on the direct sum $\widehat E=TM\oplus E$ contained in $\FRo(\widehat E)$,
where $\widehat E$ is endowed with the semi-Riemannian structure $\hat g$
of index $\bar r$ given by the orthogonal direct sum of $g$ and $g^0$.
We denote by $\nabla$ and $\overline\nabla$
respectively the Levi-Civita connections of $g$ and $\bar g$, by $\widehat\nabla$ the connection on $\widehat E$
compatible with $\hat g$ with components $\nabla$, $\nabla^{E^0}$ and $\alpha^0$ (see Remark~\ref{thm:remonly3components}),
by $\nablao$ the connection induced by $\nabla$ and $\nabla^0$ on the vector bundle $TM^*\otimes TM^*\otimes E^0$ and
by $R$ and $R^0$ the curvature tensors of $\nabla$ and $\nabla^0$ respectively.

If $(\overline M,\overline\nabla,\overline P)$ is infinitesimally homogeneous with characteristic
tensors $\overline\Tor_\oo$, $\overline R_\oo$ and $\overline{\In{}}_\oo$ then obviously $\overline\Tor_\oo=0$;
moreover, since the connection $\overline\nabla$ is compatible with the semi-Riemannian metric $\bar g$, the
linear map $\overline{\In{}}_\oo$ takes values in $\so(\bar n-\bar r,\bar r)/\mathfrak g$. Recall
that if $Z$ is any $\bar n$-dimensional real vector space endowed with a $G$-structure then we have
``versions'' $\overline R_Z$ and $\overline{\In{}}_Z$ of the characteristic tensors
$\overline R_\oo$ and $\overline{\In{}}_\oo$ of $(\overline M,\overline\nabla,\overline P)$
(recall Remark~\ref{thm:versions}).
\begin{teo}\label{thm:grotescometrico}
Assume that $(\overline M,\overline\nabla,\overline P)$
is infinitesimally homogeneous with characteristic tensors $\overline\Tor_\oo$, $\overline R_\oo$ and $\overline{\In{}}_\oo$
and that the equalities:
\begin{gather}
\label{eq:outroGauss}\begin{aligned}
\hat g_x\big(\overline R_{\widehat E_x}(v,w)u,z\big)=g_x\big(R_x(v,w)u,z\big)
&-g^0_x\big(\alpha^0_x(w,u),\alpha^0_x(v,z)\big)\\
&+g^0_x\big(\alpha^0_x(v,u),\alpha^0_x(w,z)\big),
\end{aligned}\\[5pt]
\label{eq:outroCodazzi}\begin{aligned}
\hat g_x\big(\overline R_{\widehat E_x}(v,w)u,e\big)&=
g^0_x\big((\nablao\alpha^0)_x(v,w,u),e\big)\\
&-g^0_x\big((\nablao\alpha^0)_x(w,v,u),e\big),
\end{aligned}\\[5pt]
\label{eq:outroRicci}\begin{aligned}
\hat g_x\big(\overline R_{\widehat E_x}(v,w)e,e'\big)=
g^0_x\big(R^0_x(v,w)e,e'\big)
&+g_x\big(A^0_x(e)\cdot v,A^0_x(e')^*\cdot w\big)\\
&-g_x\big(A^0_x(e)\cdot w,A^0_x(e')\cdot v\big),
\end{aligned}\\[5pt]
\label{eq:outraInnerTorsion}\overline{\In{}}_{\widehat E_x}(v)=\In{\widehat P}_x(v),
\end{gather}
for all $x\in M$, $v,w,u,z\in T_xM$, $e,e'\in E^0_x$, where $A^0$ is defined by \eqref{eq:g0alpha0A0}.
Then, for every $x_0\in M$, every $y_0\in\overline M$
and every $G$-structure preserving linear map $\sigma_0:\widehat E_{x_0}\to T_{y_0}\overline M$ there
exists a $G$-structure preserving local solution $(f,S)$ for the isometric
immersion problem with data $\nabla^0$, $\alpha^0$, $g^0$ whose domain is an open neighborhood
$U$ of $x_0$ in $M$, $f(x_0)=y_0$ and $\sigma_0=L_{x_0}$, where $L$ is as in \eqref{eq:defLfromfS}.
Moreover, if $M$ is (connected and) simply-connected and $(\overline M,\overline\nabla)$ is geodesically complete
then there exists a unique global solution $(f,S)$ for the isometric immersion problem with data
$\nabla^0$, $\alpha^0$, $g^0$ such that $f(x_0)=y_0$ and $\sigma_0=L_{x_0}$.
\end{teo}
\begin{proof}
Apply Theorem~\ref{thm:grotescao}, observing that \eqref{eq:outroGauss} implies \eqref{eq:Gausschar},
\eqref{eq:outroCodazzi} implies both \eqref{eq:Codazzichar} and \eqref{eq:Codazzi2char}, \eqref{eq:outroRicci} implies
\eqref{eq:Riccichar}, \eqref{eq:torsioneq1char} holds trivially and \eqref{eq:torsioneq2char} follows from the
symmetry of $\alpha^0$.
\end{proof}

\begin{example}[semi-Riemannian manifolds with constant sectional curvature]
Assume that the semi-Riemannian manifold $(\overline M,\bar g)$ has constant sectional curvature equal to $c\in\R$
and set $G=\Or(\bar n-\bar r,\bar r)$, $\overline P=\FRo(T\overline M)$, so that $(\overline M,\overline\nabla,\overline P)$ is
infinitesimally homogeneous (recall Example~\ref{exa:constseccurvinfhom}).
In this case, setting $\widehat P=\FRo(\widehat E)$, Theorem~\ref{thm:grotescometrico} reproduces
the classical Fundamental Theorem of Isometric Immersions (see, for instance, \cite{Dajczer,Tenenblat}). More explicitly,
equations \eqref{eq:outroGauss}, \eqref{eq:outroCodazzi} and \eqref{eq:outroRicci} become, respectively, the standard Gauss, Codazzi and Ricci equations
(observe that the lefthand side of \eqref{eq:outroCodazzi} and \eqref{eq:outroRicci} vanish); moreover, equation \eqref{eq:outraInnerTorsion}
(whose lefthand side also vanishes) says that the connection $\widehat\nabla$ is compatible with the semi-Riemannian structure $\hat g$ of $\widehat E$,
so that the equation is always satisfied.
\end{example}

\begin{example}[semi-Riemannian K\"ahler manifolds with constant holomorphic sectional curvature]
Let $\overline J$ be a $\overline\nabla$-parallel $\bar g$-antisymmetric almost complex structure on $\overline M$,
so that $(\overline M,\bar g,\overline J)$ is a semi-Riemannian K\"ahler manifold. Assume that $(\overline M,\bar g,\overline J)$ has constant holomorphic
sectional curvature equal to $c\in\R$. Setting $G=\Ur\big(\frac12\,\bar r,\frac12(\bar n-\bar r)\big)$ and
$\overline P=\FRu(T\overline M)$ then $(\overline M,\overline\nabla,\overline P)$ is infinitesimally homogeneous
(recall Example~\ref{exa:constholcurvinfhom}). Consider a $\nabla$-parallel $g$-antisymmetric almost complex structure $J$ on $M$ and
a $\nabla^0$-parallel $g^0$-antisymmetric almost complex structure $J^0$ on $E^0$; set $\widehat J=J\oplus J^0$
and $\widehat P=\FRu(\widehat E)$. In this case, Theorem~\ref{thm:grotescometrico} becomes a well-known result about
existence of isometric immersions into K\"ahler manifolds of constant holomorphic sectional curvature. Notice that
the lefthand sides of \eqref{eq:outroCodazzi} and \eqref{eq:outraInnerTorsion} vanish (but the lefthand side of
\eqref{eq:outroRicci} does not!); equation \eqref{eq:outraInnerTorsion} says that both $\hat g$ and $\widehat J$ are
$\widehat\nabla$-parallel, which happens if and only if $\alpha^0_x:T_xM\times T_xM\to E^0_x$ is complex bilinear,
for all $x\in M$.
\end{example}

Recent advances in the theory of minimal and constant mean curvature submanifolds (see for instance \cite{AbrRos, Benoit1, MeeRos})
have triggered an increasing interest in the study of immersions into product manifolds (\cite{HofLirRos, SaE}).
We discuss the case of isometric immersions into product of manifolds with constant sectional curvature.
\begin{example}[products of manifolds with constant sectional curvature]\label{exa:productimm}
Let us consider semi-Riemannian manifolds $(\overline M_1,\bar g_1)$, $(\overline M_2,\bar g_2)$ such that
$\Dim(\overline M_i)=\bar n_i$, the index of $\bar g_i$ is $\bar r_i$ and the sectional curvature of
$(\overline M_i,\bar g_i)$ is constant and equal to $c_i$, $i=1,2$. Assume that $(\overline M,\bar g)$ is the orthogonal
cartesian product of $(\overline M_1,\bar g_1)$ and $(\overline M_2,\bar g_2)$, so that $\bar n=\bar n_1+\bar n_2$
and $\bar r=\bar r_1+\bar r_2$. Let $\R^{\bar n}=\R^{\bar n_1}\oplus\R^{\bar n_2}$ be endowed with the orthogonal direct
sum of the standard Minkowski inner product of index $\bar r_1$ in $\R^{\bar n_1}$ and the standard Minkowski inner product
of index $\bar r_2$ in $\R^{\bar n_2}$. Let $G\cong\Or(\bar r_1,\bar n_1-\bar r_1)\times\Or(\bar r_2,\bar n_2-\bar r_2)$
be the Lie group of linear isometries of $\R^{\bar n}$ that preserve $F_0=\R^{\bar n_1}\oplus\{0\}$ and set
$\overline P=\FRo\big(T\overline M,(T\overline M_1)\times\overline M_2\big)$ (recall the notation introduced in Example~\ref{exa:calcInPOrEF}).
Then $(\overline M,\overline\nabla,\overline P)$ is infinitesimally homogeneous (this is obtained by applying
the product construction explained in Example~\ref{exa:prodMPs}). If $Z$ is a $\bar n$-dimensional real vector
space then a $G$-structure on $Z$ can be identified with a pair $(\langle\cdot,\cdot\rangle,Z_1)$, where
$\langle\cdot,\cdot\rangle$ is a nondegenerate symmetric bilinear form of index $\bar r$ on $Z$ and $Z_1$ is a subspace of $Z$
on which $\langle\cdot,\cdot\rangle$ is a nondegenerate symmetric bilinear form of index $\bar r_1$. The ``versions'' of the
characteristic tensors on $Z$ become $\overline{\In{}}_Z=0$ and:
\begin{multline*}
\overline R_Z(z_1+z'_1,z_2+z'_2)(z_3+z'_3)=c_1\big(\langle z_2,z_3\rangle z_1-\langle z_1,z_3\rangle z_2\big)\\
+c_2\big(\langle z'_2,z'_3\rangle z'_1-\langle z'_1,z'_3\rangle z'_2\big),
\end{multline*}
for all $z_1,z_2,z_3\in Z_1$ and all $z'_1,z'_2,z'_3\in Z_1^\perp$.
\end{example}

In Example~\ref{exa:productimm} one can consider for instance the case where $\overline M$ is an orthogonal product of a space
form with the real line $\R$, obtaining the isometric immersion theorems of \cite{Benoit1}.

\begin{example}
An isometric immersion theorem into Lie groups endowed with left-invariant semi-Riemannian metrics
is easily obtained from Example~\ref{exa:affimmgroup} by considering a linear isometry $p_0:\R^{\bar n}\to\mathbf g$
and an orthonormal global smooth frame $s:M\to\FRo(\widehat E)$.
\end{example}

\begin{rem}
A different set up for isometric immersions into Lie groups has been discussed in \cite{JorgeHerbert} where
the author considers isometric immersions into some classes of solvable and nilpotent Lie groups.
In our terminology, the result in \cite{JorgeHerbert} is associated to the $G$-structure on the target manifold determined
by an orthonormal basis of the center of the Lie algebra. Such $G$-structure is not in general infinitesimally homogeneous,
so the result in \cite{JorgeHerbert} is not a direct corollary of Theorem~\ref{thm:grotescometrico}.
\end{rem}

\begin{example}
Assume that $(\overline M,\bar g)$ is a three-dimensional Riemannian manifold as in Example~\ref{exa:Benoit}.
Let $\overline P$ denote the $G$-structure of positively oriented linear isometries that send $(1,0,0)$ to $\xi$,
where $G\cong\SO(\R^2)$ is as in Example~\ref{exa:Benoit}; then $(\overline M,\overline\nabla,\overline P)$
is infinitesimally homogeneous. Assume that $(M,g)$ is an oriented two-dimensional Riemannian manifold and that $E=M\times\R$ is the trivial
linear bundle over $M$. Consider a global smooth unitary section of $\widehat E$, which is determined
by a smooth vector field $X$ on $M$ and a smooth function $\nu:M\to\R$ with $g_x\big(X(x),X(x)\big)+\nu(x)^2=1$,
for all $x\in M$. We consider the $G$-structure $\widehat P$ on $\widehat E$ of positively oriented linear
isometries that send $(1,0,0)$ to $(X,\nu)$. In this context, a local solution $(f,S)$ for the isometric
immersion problem with data $\nabla^0$, $\alpha^0$, $g^0$ is $G$-structure preserving if
$\xi\big(f(x)\big)=\dd f_x\big(X(x)\big)+\nu(x)S_x(1)$, for all $x\in M$ (notice that $x\mapsto S_x(1)$ is a unit
normal vector field). In this context, we have (omitting the obvious zero corresponding to the covariant derivative
of the metric tensor):
\[\overline{\In{}}_{\widehat E_x}(v)=\tau v\times\big(X(x),\nu(x)\big)
=\tau\big({-\nu(x)J_x(v)},\langle J_x(v),X(x)\rangle\big),\]
for all $x\in M$, $v\in T_xM$,
where $J_x:T_xM\to T_xM$ denotes the canonical complex structure of the two-dimensional oriented Riemannian manifold
$(M,g)$. Moreover:
\[\In{\widehat P}_x(v)=\widehat\nabla_v(X,\nu)=\Big(\nabla_vX+\nu(x)A^0_x(1)\cdot v,\alpha^0_x\big(v,X(x)\big)+\dd\nu_x(v)\Big).\]
In this case, Theorem~\ref{thm:grotescometrico} reproduces the isometric immersion theorem of \cite{Benoit}.
\end{example}

\end{section}

\end{document}